\documentclass[12pt]{article}

\newcommand  {\Ebar} {{\mbox{\rm$\mbox{I}\!\mbox{E}$}}}
\newcommand  {\Rbar} {{\mbox{\rm$\mbox{I}\!\mbox{R}$}}}

\newcommand  {\Nbar} {{\mbox{\rm$\mbox{I}\!\mbox{N}$}}}

\newsavebox{\zzzbar}
\sbox{\zzzbar}
  {\setlength{\unitlength}{0.9em}
  \begin{picture}(0.6,0.7)
  \thinlines
  \put(0,0){\line(1,0){0.6}}
  \put(0,0.75){\line(1,0){0.575}}
  \multiput(0,0)(0.0125,0.025){30}{\rule{0.3pt}{0.3pt}}
  \multiput(0.2,0)(0.0125,0.025){30}{\rule{0.3pt}{0.3pt}}
  \put(0,0.75){\line(0,-1){0.15}}
  \put(0.015,0.75){\line(0,-1){0.1}}
  \put(0.03,0.75){\line(0,-1){0.075}}
  \put(0.045,0.75){\line(0,-1){0.05}}
  \put(0.05,0.75){\line(0,-1){0.025}}
  \put(0.6,0){\line(0,1){0.15}}
  \put(0.585,0){\line(0,1){0.1}}
  \put(0.57,0){\line(0,1){0.075}}
  \put(0.555,0){\line(0,1){0.05}}
  \put(0.55,0){\line(0,1){0.025}}
  \end{picture}}
\newcommand{\Zbar}{\mathord{\!{\usebox{\zzzbar}}}}

\newsavebox{\uuunit}
\sbox{\uuunit}
    {\setlength{\unitlength}{0.825em}
     \begin{picture}(0.6,0.7)
        \thinlines
        \put(0,0){\line(1,0){0.5}}
        \put(0.15,0){\line(0,1){0.7}}
        \put(0.35,0){\line(0,1){0.8}}
       \multiput(0.3,0.8)(-0.04,-0.02){12}{\rule{0.5pt}{0.5pt}}
     \end {picture}}

\newcommand{\QED}{{\hspace*{\fill}\rule{2mm}{2mm}\linebreak}}

\newtheorem{lemma}{Lemma}[section]
\newtheorem{proposition}{Proposition}[section]
\newtheorem{theorem}{Theorem}[section]
\newtheorem{definition}{Definition}[section]

\newcommand{\Z}{\Zbar}
\newcommand{\R}{\Rbar}
\newcommand{\N}{\Nbar}

\newcommand{\E}{\Ebar}

\begin{document}

\setlength{\textheight}{21cm}
\title{The Restriction of the Ising Model to a Layer
 }
\author{C. Maes\thanks{Onderzoeksleider FWO, Flanders. \
                        Email: christian.maes@fys.kuleuven.ac.be }  \\
    F. Redig\thanks{Postdoctoraal onderzoeker FWO, Flanders. \ 
                        Email: frank.redig@fys.kuleuven.ac.be} \\
    A. Van Moffaert\thanks{Aspirant FWO, Flanders. \
                       Email: annelies.vanmoffaert@fys.kuleuven.ac.be }
                 \thanks{Address :
    Celestijnenlaan 200D, B-3001 Leuven, Belgium } \\
    Instituut voor Theoretische Fysica  \\
         K.U.Leuven}
\maketitle

\begin{abstract}
We discuss the status of recent Gibbsian descriptions of the restriction
(projection) of the Ising phases to a layer.
We concentrate on the projection of the two-dimensional low temperature Ising
phases for which we prove a variational principle.
\end{abstract}
\vspace{3mm}



\section{Introduction}

In this paper we study the restriction of the two-dimensional Ising model to a
(one-dimensional) layer. The restriction of the plus (or minus) phase is known
to be non-Gibbsian below the critical temperature, see \cite{Sch}, \cite{FP}.
Following suggestions of Dobrushin it was recently shown that this restriction
is in fact weakly Gibbsian, see \cite{D}, \cite{DS}, \cite{MV}, \cite{MRV}.
We state and discuss various recent results on these restrictions. Using
elementary methods, we rederive results on weak 
Gibbsianness below the critical
temperature and the results on Gibbsianness above the critical temperature
or in a magnetic field. We also add the result that further random decimations
of the weakly Gibbsian restriction are again Gibbsian.
Finally, we prove the existence of thermodynamic functions (energy and free
energy density) and we discuss a variational principle for the weakly Gibbsian
measure.

The study of restrictions of Gibbs
measures (on $d+1$-dimensional configurations) to a sublayer (of dimension $d$)
can be motivated in various ways.
First of all they are interesting test
cases for an extended Gibbsian description of non-Gibbsian states.
Since about ten years various examples of non-Gibbsian states have been
produced.  Some of these go back to the work of Griffiths, Pearce and
Israel, \cite{GP1}, \cite{Isr}, and have become (in)famous as so called
renormalization group pathologies, \cite{EFS}.
Dobrushin's program
tries to understand the non-Gibbsianness as coming from a perhaps
too strict requirement on the potential.  If, like is the case for
unbounded spins, one asks for a potential which is summable on what
are typical configurations for the state, one can get at least some
effective interaction or physically relevant parametrization of e.g. the
images under transformations of states.


A second motivation comes from the theory of interacting particle systems. One
of the questions is to see under what conditions the invariant measures of a dynamics
are Gibbsian. The simplest scenario is found for so called probabilistic
cellular automata (PCA). These are stochastic dynamics for lattice spin systems
under which the spins are updated synchronously in discrete time. If one starts
a PCA (with positive transition probabilities) from an invariant measure, 
then the distribution of space-time configurations
(configurations on the space-time lattice) turns out to be a Gibbs measure
see \cite{LMS}. Therefore, the invariant measure itself is a restriction of a
$d+1$-dimensional
Gibbs measure 
to a $d$-dimensional hypersurface.

In fact, the question in \cite{Sch} about the Ising restrictions
came quite
naturally after it was found in \cite{LMS} that for high noise dynamics 
the unique stationary state is Gibbsian as
follows from considering the restriction of a Gibbs state in the regime of
complete analyticity.  The question about the Gibbsian nature of stationary
states has of course been considered before, see e.g. \cite{Kun},
\cite{Kun1}, \cite{MV93}, \cite{MV94} and references therein.  Understanding
the locality of the time reversal with respect to the stationay state plays
a crucial role in these.  
While
most of these problems are  still open (for general PCA's), we feel that the
Dobrushin program gives new inspiration towards a (weakly) Gibbsian description 
of these
invariant measures if considered as restrictions of Gibbs measures.

A third motivation can be found in the study of surfaces, models on 
half-planes
with random boundary conditions (wetting phenomena). The restriction of the two
dimensional Ising model to a layer can of course be seen as a surface state 
with respect
to the two dimensional Ising measure. The problem of finding an interaction for
this restriction consists in finding the interaction between spins at the 
boundary of
an Ising sample (as a function of the configuration at this boundary). 
The relation with wetting is of a more technical
nature. It turns out that in the study of the restriction of Ising measures, 
one is
quickly confronted with questions like how far the influence of a configuration
on a sublayer is felt in the bulk of the system. The wetting-context can be 
used to get a
very useful intuitive picture of why, at low temperatures, the restrictions 
of the 
plus phase of the Ising model are not Gibbsian. This was made precise in 
\cite{EFS},
see also \cite{FP}. An interesting further question (related to the convergence
properties of the potential for the restriction) is to see whether there is good
decay of correlations close to the surface when on this surface we impose a
'typical surface-configuration' (i.e. a sample of the restriction). 
We will answer this
question in Section 4 (Proposition \ref{prop1}). 

The paper is organized as follows: in Section 2 we introduce basic notations
and definitions. We introduce the 'telescoping' potential 
(\`{a} la Kozlov, \cite{Koz}) in Section 3 and
discuss its summability properties. In Section 4 we give an overview of the
results on the restrictions of the Ising model. In Section 5 we prove the results
of Section 4 using the telescoping potential, and finally in Section 6 we discuss
the variational principle.

\section{Definitions and Notations}
\setcounter{equation}{0}
\subsection{Configuration space}

We consider
the regular $r-$dimensional lattice
$\Z^r$ and 
denote by 
$\mathcal{L}:= \{V, |V| < \infty\}$  
the set of finite subsets of $\Z^r$.  The complement of a set
$V\subset \Z^r$ is $V^c = \Z^r \setminus V$.
For two sites $x,y \in \Z^r$ we define
\begin{equation}
\label{dist}
|x-y| := \sum_{i=1}^r |x_i - y_i|.
\end{equation}

The state space is $\Omega := \{+1,-1\}^{\Z^r}$ 
and
its elements (= configurations) are denoted by
greek letters
$\eta,\omega,\sigma,\xi,\ldots$. 
The value of $\omega $ at a site $i \in \Z^r$ is written as
$\omega(i) $.
For $V \in \mathcal{L}$, $\sigma \in \Omega$ we define

\begin{equation}
\sigma^V (x)   = 
 \left\{
     \begin{array}{ll}
      \sigma (x)  & \mbox{if $ x \in V$} \\
	 +1       & \mbox{if $ x \not\in V $}.
     \end{array}
 \right. 
\end{equation} 
For $\sigma, \eta \in \Omega$ we define $\sigma_V\eta_{V^c}$
to be the configuration
\begin{equation}
\sigma_V \eta_{V^c} (x)  = 
 \left\{
   \begin{array}{l}
    \sigma(x) \quad x \in V \\
    \eta(x)   \quad x \not \in V.
   \end{array}
 \right. 
\end{equation}
The restriction of $\Omega$ to a volume $V \subset \Z^r$ is denoted by
$\Omega_V := \{+1,-1\}^V$ and we write
$\sigma_V \in \Omega_V$ for the restriction of a configuration
$\sigma$.

On $\Omega$ we have the natural action of translations $\tau_a$,
$a \in \Z^r$ defined
by $ \tau_a \eta(i) := \eta(i-a), i \in \Z^r$.
The $\sigma$-algebra generated by the evaluation maps
$X_i, X_i(\omega) :=  \omega(i), i \in M$  is written as
$ \mathcal{F}_M = \sigma\{ X_i, \ i\in M\}$.  
When $M = \Z^r$, we set $\mathcal{F} := \mathcal{F}_{\Z^r}$.
The {\it tail field} $\sigma$-algebra $\mathcal{T}^{\infty}$ is defined as
\begin{equation}\label{tail}
\mathcal{T}^{\infty} := \cap_{V \in \mathcal{L}}
 \mathcal{F}_{V^c}.
\end{equation}
The configuration space $\Omega$ is a compact metric space in the product
topology.
A function $f$ on $\Omega $ is called local if it depends only on a finite 
number
of coordinates, i.e. there is a $V \in \mathcal{L}$ such that
$f(\eta) = f(\zeta)$ whenever $\eta_V = \zeta_V$.
The minimal set $V$ such that this holds is called the 
{\em dependence set} 
of the function.

\begin{definition}
\label{defrightcont}
A function $f : \Omega \rightarrow \R$
is called {\em right-continuous} in $\sigma \in \Omega$ if
 \begin{equation}
 \label{cont}
 f(\sigma) = 
 {\rm lim}_{V \uparrow \Z^r} 
    f(\sigma^V) .
 \end{equation}
\end{definition}

On $\Omega$ we have the pointwise order $\eta \preceq \omega$ if $\eta(x) \leq
\omega(x)$ $\forall x \in \Z^r$.
A function $f: \Omega \rightarrow \R$
is called {\em monotone non-decreasing} if
for all $\eta, \xi \in \Omega$, $ \eta \preceq \xi$ implies 
$f(\eta) \leq f(\xi)$.

\subsection{Potentials and specifications}

\begin{definition}[cf \cite{FP}]
\label{specification}
A  local specification $\Gamma$ on $\mathcal{L}$ is a family of probability
kernels $\Gamma = \{ \gamma_{V},  V \in \mathcal{L} \}$
on $(\Omega, \mathcal{F})$, such that the following hold:
\begin{enumerate}
\item
 $\gamma_{V}(\cdot|\omega)$ is a probability measure on $(\Omega,
 \mathcal{F})$ for all $\omega \in \Omega$;
\item
 $\gamma_{V}(F|\cdot)$ is $\mathcal{F}_{V^c}$ -measurable for all
 $F \in \mathcal{F}$;
\item
 $\gamma_{V}(F|\omega) = 1_F (\omega)$ if
 $F \in \mathcal{F}_{V^c}$;
\item
 $\gamma_{V_2} \gamma_{V_1} = \gamma_{V_2} $ if
 $V_1 \subset V_2$.
\end{enumerate}
\end{definition}

\begin{definition}
\label{meas-spec}
A probability measure $\mu$ is consistent with a specification $\Gamma$ (or vice
versa), notation $\mu \in \mathcal{G}(\Gamma)$,
if $\forall V \in \mathcal{L}$
\begin{equation}
\mu \ \gamma_V = \mu.
\end{equation}
\end{definition}

A specification $\Gamma$ is said to be {\it translation invariant} if
$ \forall a \in \Z^r$, $\forall V \in \mathcal{L}$, $\forall \omega \in \Omega$
and for all bounded measurable functions $f$
\begin{equation}
\gamma_V(f \circ \tau_a|\omega) = \gamma_{V + a} (f|\tau_a \omega).
\end{equation}
It is customary to slightly abuse the notation (and to circumvent property
3. above) by writing $\gamma_V(\sigma|\omega) = 
\gamma_V(\sigma_V|\omega_{V^c})$
if one means to take configurations $\sigma$ and $\omega$ identical on 
$V^c$. 
One should think of 
 $\gamma_V(\sigma|\omega)$,  as
the probability to find $\sigma$ in $V$ given $\omega$ outside of
$V$.

Property 4 of Definition (\ref{specification}) is called
self-consistency and is most important in characterising equilibrium.
It suggests constructing probability measures 
$\nu\in \mathcal{G} (\Gamma )$
as weak limits of $\gamma_{V}(\cdot|\omega),$ some $\omega \in \Omega,
V \uparrow \Z^r$ (perhaps along a subsequence). Such weak limits
automatically exist by compactness but their consistency with the specification
is only immediate if
$\gamma_V(f|\cdot)$ is a continuous function for
all continuous $f$.  One then deals with a so-called Feller (or quasilocal)
specification.  This is not the context of the Dobrushin program
where more general specifications have to be considered and hence
that $\mathcal{G}(\Gamma) \not= \emptyset$ is not obvious in general. 
 
In a Gibbsian formalism  
one considers a special class of specifications, the
so called Gibbsian specifications which are of the 
Boltzmann-Gibbs form:
\begin{equation}
\gamma_V(\sigma|\omega) = \frac{1}{Z_V(\omega)}
 \exp \{ - \sum_{A\cap V \neq \emptyset} 
 U(A, \sigma_V \omega_{V^c}) \}
\end{equation}
where 
$Z_V(.) \in\mathcal{F}_{V^c}$ 
is a normalization factor
\begin{equation}
Z_V(\omega) = \sum_{\sigma_V \in \Omega_V} 
\exp \{- \sum_{A\cap V \neq
\emptyset} U(A, \sigma_V\omega_{V^c})\}
\end{equation}
and $U(A,.)$ is an ``interaction potential'':

\begin{definition}
A potential $U$ is a real-valued function on $\mathcal{L} \times \Omega$
\begin{equation}
U: \ \mathcal{L} \times \Omega  \rightarrow \R
\end{equation}
such that $U(A,\cdot)$ $\in$ $\mathcal{F}_A$ for all $A
\in \mathcal{L}$ (put $U(\emptyset,\cdot) = 0$).
\end{definition}

A potential $U$ is {\it translation invariant} if
$ \forall A \in \mathcal{L}$, $a \in \Z^rS$, $\eta \in \Omega$
\begin{equation}
\label{Z}
U(A,\eta) = U(A +a,\tau_a \eta).
\end{equation}

\begin{definition}
\label{potential}
\begin{enumerate}
\item
A potential $U$ is {\em convergent} in $\omega\in \Omega$ 
if for all $V \in
\mathcal{L} $
\begin{equation}
\label{conv}
\sum_{A \cap V \neq \emptyset}
  U(A,\omega) 
\end{equation}  
is well-defined.
We always understand an infinite sum $\sum_A a_A$ as well-defined
when $\exists a< \infty$  $\forall \epsilon >0$ 
$\exists V_0 \in \mathcal{L}$
so that $\forall V \in \mathcal{L}$, $V \supset V_0$
\[
|\sum_{A\subset V} a_A  - a|  \leq \epsilon.
\]

\item
A potential $U$ is {\em absolutely convergent} in $\omega \in \Omega$ if
for all $V \in \mathcal{L}$
\begin{equation}
\label{absconv}
 \sum_{A \cap V \neq \emptyset}
   |U(A,\omega)|  < \infty .
\end{equation}

\item
A potential $U$ is {\em uniformly absolutely convergent} if
for all $V \in \mathcal{L}$
\begin{equation}
\label{uniabs}
 \sum_{A \cap V \neq \emptyset}
   {\rm sup}_{\omega \in \Omega}|U(A,\omega)|  < \infty. 
\end{equation}
\end{enumerate}  
\end{definition}

Let $U$ be a potential and suppose that there exists  a  set 
$\bar{\Omega}_U $ in the tail field of points of convergence of $U$ 
($\bar{\Omega}_U \in \mathcal{T}^\infty$
and $\forall \omega \in \bar{\Omega}_U$ and $\forall V \in \mathcal{L}$
the sum $\sum_{A \cap V \not= \emptyset} U(A,\omega)$
is well-defined).
Then, for every 
 $V \in \mathcal{L}$
and every $\omega \in \bar{\Omega}_U$
we can introduce the finite volume Gibbs measure
\begin{equation}
\label{finiteGibbs}
\mu_V^{\omega,U} (\xi)  = 
 \left\{
   \begin{array}{ll}
     \frac{1}{Z_V(\omega)} \exp \{-\sum_{A\cap V \neq \emptyset}
      U(A, \xi_V \omega_{V^c}) \} & \mbox{if} \ 
        \xi = \xi_V \omega_{V^c}  \\
      0  & \mbox{otherwise} .
   \end{array}
 \right.  
\end{equation}
( We ask for $\bar{\Omega}_U$ to be in the tail field to make sure that 
$Z_V(\omega)$ is well-defined.)
Factors of temperature or {\it a priori} weights (reference measure) are
supposed to be contained in the potential.
The Dobrushin operator is then defined by taking expectations with respect to
(\ref{finiteGibbs}):

\begin{equation}
\label{Doop}
R^U_V(f)(\omega) := 
\int f(\xi) \mu_V^{\omega,U}(d\xi)
\end{equation}
mapping bounded measurable functions $f$ on $\Omega$  
to functions $R^U_V(f)$ on $\bar{\Omega}_U$.


%

\begin{definition}
\label{defweakly}
A probability measure $\mu$ on $(\Omega,\mathcal{F})$ is {\em weakly
Gibbsian}
if there exists a potential $U$ and a tail field set
$\Omega_U $ of points of absolute convergence of $U$ (cf. (\ref{absconv}))
such that
\begin{enumerate}
\item
 $\mu(\Omega_U) = 1$ ;
\item
 $\forall V \in \mathcal{L}$, $\forall \mathcal{B}
 \in\mathcal{F}_{V^c}$
 and for every
 bounded measurable function $f$, 
\begin{equation}
  \int_{\mathcal{B}} f \ d\mu = \int_{\mathcal{B}} 
  R^U_V(f)
  d\mu .
\end{equation}
\end{enumerate}
\end{definition}
A somewhat less stringent definition of weak Gibbsianness is obtained by asking that
there is a tail field set $\bar{\Omega}_U$ of points of  convergence of
$U$ such that 1. and 2. of Definition \ref{defweakly} hold.

If the potential in Definition (\ref{defweakly}) is uniformly absolutely
convergent, then $\nu$ is a ({\it bona fide}) Gibbs measure.

\section{Vacuum and Telescoping Potential.}
\setcounter{equation}{0}
\label{section3}
In this section we shall introduce the so-called telescoping potential
which will be a very useful tool in the study of the restrictions
of the Ising model. 
A natural potential associated to a specification is the so-called
vacuum potential (in our case the vacuum will always be the configuration
of all plusses). To construct this potential, start from
\begin{equation} \label{sta}
H_V(\xi) := \ln
\frac{\gamma_V(+|+)}{\gamma_V(\xi^V|+)}
\end{equation}
and write
\begin{equation}\label{vache}
H_V (\xi ) = \sum_{A\subset V } v(A,\xi ).
\end{equation}
This last formula can be inverted (M$\mbox{\o}$bius formula) and we get:
\begin{equation}\label{vac}
v(A,\xi) = \sum_{V\subset A} (-1)^{|A\setminus V|} H_V(\xi).
\end{equation}
The inversion is rapidly checked by remembering that
$\sum_{V\subset R \subset
B} (-1)^{|R|} = (-1)^{|V|} \delta_{B,V}$ where 
$\delta_{\cdot}$ is the Kronecker
delta. 
This potential $v$ is called {\it vacuum} because it has the property that
$v(A,\omega) =0$ whenever $\omega_i=+1$ for some $i\in A$ 
and $\xi(i) = + $.  This
follows easily from (\ref{vac}) using that $H_V(\xi) = H_{V\setminus
i}(\xi)$ if $i \in V \subset A$ and $\xi(i) =+$.  It is straightforward to check that the
vacuum potential
is the unique potential having this property and that it is 
translation invariant if $\Gamma$ is.  
In \cite{MRV} it is proved that the
vacuum potential $v$ is convergent and consistent with the specification 
$\Gamma$, i.e. $\forall V \in \mathcal{L}$ $\forall \omega, \xi \in
\Omega$
\begin{equation}
\label{vac-gamma}
\gamma_V(\xi|\omega)=\frac{1}{Z_V(\omega)} 
 \exp \left[ - \sum_{A\cap V \not= \emptyset} 
 v(A,\xi_\Lambda\omega_{\Lambda^c}) \right]
\end{equation}
if and only if the specification is right-continuous.

A possible problem with this vacuum potential is that it may not be
{\it absolutely} convergent
(even when the specification is right-continuous). Therefore, in
order to obtain absolute convergence
Kozlov, in \cite{Koz}, introduces another kind of potential. We now
present
a simplified version of it which, for the occasion, we would like to call a
telescoping potential.
For this we turn again to (\ref{sta}) and we write
\begin{equation}  \label{beg}
\exp(-H_V(\xi)) =
\frac{\gamma_V(\xi^V|+)}{\gamma_V(+|+)} = \prod_{s=1}^n
\exp F_s(\xi_{i_1}\cdots \xi_{i_s})
\end{equation}
where we have lexicographically ordered the $n=|V|$ sites in $V$ according
to $i_1 < i_2 < \cdots < i_n$ and
\begin{equation}\label{iff}
F_s(\xi_{i_1}\cdots \xi_{i_s}) := \ln \frac{\gamma_V(\xi_{i_1}\cdots
\xi_{i_s} + \cdots +|+)}{\gamma_V(\xi_{i_1}\cdots \xi_{i_{s-1}} + \cdots
+|+)}. \end{equation}
We can now order the sites $i_1, \ldots, i_s = \{j\leq i, j\in V\}$
according to their
distance from the `largest' site $i_s=i$.  For this purpose we consider
for every $i\in V$ the sequence of increasing volumes
$L_{i,m}$ with $L_{i,m} := \{j \in \Z^r : j \leq i, |j-i| \leq m\},
m=0, 1, \ldots$.  We thus have the partition
\begin{equation} \label{partit}
\{j\leq i, j\in V\} = \cup_{m=1}^{m(i,V)} V_{i,m} \setminus V_{i,m-1} \cup
V_{i,0}
\end{equation}
with $V_{i,m} := L_{i,m} \cap V$ and $m(i,V) \equiv \max_{j\leq i, j\in
V} |i-j|$.   Correspondingly, $F_s(\xi_{i_1}\cdots \xi_{i_s}) =
F_i(\xi^{V_{i,m(i,V)}})
$ can be
further telescoped as
\begin{equation}   \label{tele}
= - \sum_{m=0}^{m(i,V)} U_{L_{i,m}}(\xi^V),
\end{equation}
with
\begin{equation} \label{ja}
U_{L_{i,m}}(\xi) :=  \ln \frac{\gamma_{L_{i,m}}(\xi^{L_{i,m-1}}|+)
\gamma_{L_{i,m}}(\xi^{L_{i,m}\setminus i
}|+)}{\gamma_{L_{i,m}}(\xi^{L_{i,m-1}\setminus i}|
+ ) \gamma_{L_{i,m}}(\xi^{L_{i,m}}|+)},
\end{equation}
for $m > 0$ and
\begin{equation} \label{janul}
U_{L_{i,0}}(\xi) := -F_i(\xi^{L_{i,0}}) = \ln \frac{\gamma_{i}(+|+)}
{\gamma_{i}(\xi^i|+)}.
\end{equation}
(Observe that $m(i,V)=0$ when
$i=i_1$ is the `first' site in $V$).
We thus define the (telescoping) potential
\begin{equation}
U(A,\xi) :=  U_{L_{i,m}}(\xi) \mbox{ see (\ref{ja}), if }
A = L_{i,m} \mbox{ for some } i\in \Z^r, m \geq 0 ;
\end{equation}
and $U(A,\xi) \equiv 0$ otherwise.

To get more insight in the potential it is
instructive
to rewrite it for the one-dimensional case. For $r=1$ the potential
$U(B,\xi)$ is non-vanishing iff $B=L_{i,m}=[i-m,i] \cap \Z, i\in \Z, m=0,
1, \ldots$ is an interval. For such a $B$ we rewrite (\ref{ja}) as
\begin{equation} \label{ja1}
U_{[j,i]}(\xi) =  \ln \frac{\gamma_{[j,i]}(\xi^{]j,i]}|+)
\gamma_{[j,i]}(\xi^{[j,i[}
|+)}{\gamma_{[j,i]}(\xi^{]j,i[}|
+ ) \gamma_{[j,i]}(\xi^{[j,i]}|+)},
\end{equation}
where we abbreviated e.g. $]j,i] \equiv \{j+1, \ldots, i\}$ for $j<i$ in
$\Z$.

Some properties of this potential are immediate.  For example, $U_{L_{i,m}}
= 0$ whenever $\xi_i=+$ or when $\xi=+$ on the set $L_{i,m}\setminus
L_{i,m-1}$. As a consequence and following (\ref{beg})-(\ref{janul}), the
Hamiltonian (\ref{sta}) is telescoped as
\begin{equation} \label{telcon}
H_V(\xi) = \sum_{A\cap V \neq \emptyset} U(A,\xi^V).
\end{equation}
Moreover, the
potential
is explicitly translation-invariant if the specification is.

From now on we will assume that the specification $\Gamma$ is 
 right-continuous.
In the sequel this specification will always be the
 monotone right-continuous specification (introduced in
\cite{FP}) consistent with the
restriction of the Ising model.

In order to verify the consistency of this telescoping potential with
the right-continuous specification, we can rewrite $U$ as a
resummation of the vacuum potential (see also \cite{Koz}). More precisely
\begin{equation}\label{verb}
U_{L_{i,m}}(\xi) = \sum_{R\ni i, R \not\subset L_{i,m-1}, R \subset
L_{i,m}} v(R,\xi)
\end{equation}
for $m > 0$, and for $m=0$ we have
$U_i(\xi) = v(i,\xi)$.

From this, one can check that $\forall V \in \mathcal{L}$ 
\begin{equation}
\label{tele-vac}
\sum_{A\cap V \not= \emptyset} U(A,\xi^V) = 
\sum_{A\cap V \not= \emptyset} v(A,\xi^V)
\end{equation}
if for every $V \in \mathcal{L}$, $U$ is absolutely convergent in 
$\xi^V$
(see \cite{Koz}).
If there exists a set $\Omega_U$ in the tail field such that
$U$ is absolutely convergent in every point $\omega$ of $\Omega_U$
then (\ref{tele-vac}) 
together with the right-continuity of $\Gamma$ and the
consistency of  $v$ with $\Gamma$ 
(see \ref{vac-gamma}) give that
\begin{equation}
\gamma_V(\xi|\omega)= \frac{1}{Z_V(\omega)} 
 \exp\left[ - \sum_{A\cap V \not= \emptyset} 
 U(A,\xi_V\omega_{V^c}) \right],
  \omega \in \Omega_U.
\end{equation}

The representation (\ref{verb}) of the telescoping potential in terms
of the vacuum potential is very useful, because one has a certain freedom in
the choice of the sets $L_{i,m}$. The only constraint on these sets
is that the obtained potential (from (\ref{verb})) is still consistent with
the specification, i.e.
\begin{equation}\label{constraint}
\forall R\in\mathcal{L}\quad\exists ! (i,m)\mbox{ such that }
R\ni i, R \not\subset L_{i,m-1}, R \subset
L_{i,m}.
\end{equation}
Indeed, when the constraint (\ref{constraint}) is satisfied we have
\begin{equation}
\sum_{A\cap V \neq \emptyset} U(A,\xi^V)
=\sum_{L_{i,m}\cap V \neq \emptyset}
\sum_{R\ni i, R \not\subset L_{i,m-1}, R \subset
L_{i,m}} v(R,\xi^V)=\sum_{A\cap V \neq \emptyset} v(A,\xi^V)
\end{equation}
and from this, together with the right-continuity of the specification, we
conclude that the potential $U$ is consistent with the specification.
The constraint (\ref{constraint}) on the sets $L_{i,m}$ is of a geometric
nature. In dimension $r=1$ we can choose e.g. $L_{i,m} = [i-g(m),i]$, where
$g(m)$ is some strictly increasing function. The freedom in the choice
of $g(m)$ permits to 'tune' a bit the convergence properties of
the potential $U$ (for our study of the restrictions it will be okay
to choose $g(m)=m$). 
In $r\geq 2$ the sets $L_{i,m}$ introduced before also satisfy this
constraint, whereas e.g. 
$L_{i,m}:= \{ j\in \Z^r : |j-i|\le m\}$ do not satisfy the constraint.

In applications, the right-continuous specification $\Gamma$
is often constructed starting from a probability measure $\nu$ such that 
$\nu \in \mathcal{G}(\Gamma)$ (cf. Definition \ref{meas-spec}).
We then know that $\nu$ is weakly Gibbsian (cf. Definition  \ref{defweakly})
if there exists a tail field set $\Omega_U $ of points of absolute convergence
of $U$ 
such that $\nu(\Omega_U) =1$.
I.e. proving that $\nu$  is weakly Gibbsian boils down to showing that

\begin{equation} \label{boil}
\sum_{j\geq i} \sum_{m \geq |i-j|} |U(L_{j,m},\xi)| < \infty, i \in \Z^r,
\end{equation}
for a full-measure (tail-)set of $\xi$'s.

\begin{proposition} \label{prop2}

Let $\Gamma$ be a local right-continuous specification, $\nu \in
\mathcal{G}(\Gamma)$ and $U$ the telescoping potential defined by
(\ref{ja}) and (\ref{janul}).
Suppose that $\exists C_1,C_2,M < \infty$and $\exists \lambda >0$ such that 
$\forall m \geq M$, $ \forall \xi \in \Omega$ $\forall j \in \Z^r$

\begin{equation} \label{hope}
|U(L_{j,m},\xi)| \leq C_1 m^{r-1} I[m\leq \ell(j,\xi)] + C_2 I[m >
\ell(j,\xi)] m^{r-1} \exp[-\lambda m],
\end{equation}
with $\ell(j,\xi) \geq 0$ possibly unbounded.

Suppose further that 
$\exists$ a translation invariant
tail field set $K, \nu(K)= 1$ 
so that $\forall \xi \in K$ $\forall j\in \Z^r$,
$\ell(j,\xi) < \infty$ and

\begin{equation} \label{tata}
|\{j \geq i : \ell(j,\xi) \geq |j-i|\}| < \infty, i\in \Z^r.
\end{equation}
Then the telescoping potential $U$ is absolutely convergent for $\xi\in K$
and $\nu$ is weakly Gibbsian.
\end{proposition}
Of course,
the left hand side of (\ref{hope})  is a local function for fixed $m$ while
the right hand side can be highly non-local as a function of $\xi$
and deals with  the dependence of the potential on $\xi$ as $m$ grows.

{\bf Remark 1 :} 
Conditions (\ref{hope}) and (\ref{tata}) may seem weird or ad hoc. 
In Section 5 we show that they are satisfied for the
one-dimensional
restriction of the plus-phase of the two-dimensional Ising model. 
In fact if $\nu$ is the
restriction to a hyperplane of a measure $\mu$
then
$U_{L_{j,m}}(\xi)$ measures the correlations between spins at sites $(j,1)$
and $(k,1), k\in L_{j,m} \setminus L_{j,m-1}$ in a measure 
$\mu^{\xi^{L_{j,m}}}$
where $\mu^{\xi^{L_{j,m}}}$ is a constrained measure obtained from 
$\mu$ and $\xi$ plays the role of a boundary condition. (This will become more
clear in Section 4).

One should think of the $\ell(j,\xi)$ as the radius of a
ball around site $j$  outside which the spins at sites $(k,1), k\leq j$ are
only
weakly correlated to the spin at site $(j,1)$ in the constrained measure
$\mu_\beta^\xi$ of (\ref{murand}).\newline
To satisfy (\ref{tata}) on a
set of full measure
for $\nu$ it suffices that $\nu[\ell(j,\xi) > u] \leq  \exp[-c u]$ for some
$c>0$.

{\bf Remark 2 :} For $r=1$, it is convenient to use the left/right symmetry
of the sets $L_{j,m}$ which are just lattice intervals $[i,j]$.  We then
ask (\ref{hope}) (which looks to the left) together with the existence of
finite $\ell^+(i,\xi)$ for which
\begin{equation} \label{hop1}
|U([i,j],\xi)| \leq C_1 I[j-i\leq \ell^+(i,\xi)] + C_2 I[j-i >
\ell^+(i,\xi)]  \exp[-\lambda |j-i|],
\end{equation}
(looking to the right). The
assumption
(\ref{tata}) in Proposition \ref{prop2} can be replaced by the requirement
that for each $i\in \Z, \xi \in K$, there
are finite $\ell^+(i,\xi), \ell^-(i,\xi) \equiv \ell(i,\xi) < \infty$ so
that for all $j > i \in \Z$, (\ref{hope}) and (\ref{hop1}) hold.
The idea is that $\ell^-(i,\xi)$
looks in the configuration $\xi$ to the left of $i$ while
$\ell^+(i,\xi)$ looks to the right.

{\bf Proof of Proposition \ref{prop2}} : We have to check (\ref{boil}).
Inserting (\ref{hope}) there are two sums to control.  The sum for $m >
\ell(j,\xi)$ is easily taken care of using the exponential decay.  The sum
over $m\leq \ell(j,\xi)$ has only a contribution if $|j-i| \leq
\ell(j,\xi)$.  We thus get that (\ref{boil}) is bounded by
\begin{equation}
C_1 \sum_{j\geq i, |j-i| \leq \ell(j,\xi)} \ell(j,\xi)^r + C_3 \sum_{j\geq
i} \exp[-\lambda |j-i|].
\end{equation}
Using assumption (\ref{tata}) this is finite and the conclusion follows
from the remarks above.
\QED

\section{The restricted Ising model}
\setcounter{equation}{0}
\subsection{The model}
\label{Ising}

We consider the standard ferromagnetic nearest neighbor Ising model on the
regular $d+1$- dimensional lattice $\Z^{d+1}, d\geq 1$. 
The symbols $\Lambda, \Lambda_n,...$ will be reserved to
indicate
finite subsets of $\Z^{d+1}$.  Their complement is $\Lambda^c = \Z^{d+1}
\setminus \Lambda$ etc.  The configuration space for the Ising model is
$\Omega' = \{+1,-1\}^{\Z^{d+1}}$
Fix $\beta \geq
0, h > 0$.  For a finite box $\Lambda\subset \Z^{d+1}$ with free (or empty)
boundary conditions, the Gibbs state $\mu_{\Lambda,\beta}^{h}$ for the
Ising model assigns a probability
\begin{equation} \label{fin}
\mu_{\Lambda,\beta}^{h}(\sigma_x, x\in \Lambda) = \frac
1{Z_{\Lambda,\beta}^{h}}
\exp[\beta\sum_{\langle xy \rangle \subset \Lambda} (\sigma_x \sigma_y -1) + h
\sum_{x\in \Lambda} \sigma_x]
\end{equation}
to an Ising spin configuration $\sigma_x = \pm 1, x\in \Lambda$. The
normalization $Z_{\Lambda,\beta}^{h}$ is the partition function (for free
boundary conditions).
The parameter $\beta$ is (proportional to) the inverse temperature and $h$
is called the magnetic field.  The first sum
in (\ref{fin}) is over the nearest neighbor pairs $\langle xy\rangle$ in
$\Lambda$.  Each site $x\in \Z^{d+1}$ has $2(d+1)$ nearest neighbors and we
write $y\sim x$ if the site $y\in \Z^{d+1}$ is a nearest neighbor of $x$.
The infinite volume Ising state $\mu_{\beta}^{h} = \lim_\Lambda
\mu_{\Lambda,\beta}^{h}$ is obtained in the thermodynamic limit as
$\Lambda\uparrow \Z^{d+1}$ along a sequence of sufficiently regular
volumes.
We can take the (weak) limit $h\downarrow 0$ of $\mu_{\Lambda,\beta}^{h}$
and
we write $\mu_\beta^+ = \mu_{\beta}^{0^+}$ for this limit. In the same
way, starting with $h < 0$ and letting $h \uparrow 0$ we can define
$\mu_\beta^-$.   We refer to all of these as Ising states.  When making no
distinction between them, we denote these states by the common symbol
$\mu_\beta$ and the corresponding random field is denoted $X=(X(x),x\in
\Z^{d+1})$.  They
are translation invariant probability measures on $(\Omega',{\cal
F}_{\Z^{d+1}})$ and they all satisfy the
Dobrushin-Lanford-Ruelle equation
\begin{equation} \label{DLR}
\mu_\beta[X(x)=\sigma_x|X(y), y\in \Z^{d+1}\setminus \{x\}](\sigma)
= \frac 1{1 + \exp(-2\beta \sigma_x\sum_{y\sim x}\sigma_y - 2h\sigma_x)}
\end{equation}
$\mu_\beta$ almost surely.   For $h=0$ and for sufficiently large
$\beta$ there are other solutions to (\ref{DLR}) (even
non-translation invariant ones  if $d\geq 2$) but we will restrict us in
what follows to the Ising states introduced above.
In particular, there is a critical value $0 < \beta_c <
\infty$ for which $\mu_\beta^+ \neq \mu_\beta^-$ whenever $\beta >
\beta_c$.
The standard methods, results and more details about the Ising model can be
found in almost any textbook on  statistical mechanics, see e.g.
\cite{Sin}, \cite{Geo} or \cite{Sim}.

We now fix a hyperplane or layer
\begin{equation} \label{plane}
\Xi = \{x=(x_1,\ldots,x_d,x_{d+1}) \in \Z^{d+1} : x_{d+1}  = 0\}
\end{equation}
which we can identify with $\Z^{d}$.  The sites in $\Xi$ are denoted by
$i, j, k, \ldots$, which, though treated as elements of $\Z^{d}$, ought to
be identified with $(i,0), (j,0), (k,0)$, $\ldots$ when considering $\Xi
\subset
\Z^{d+1}$. Finite subsets of $\Xi$ are written as $V, V_n, A \ldots$ and
by $V^c$ we mean the complement of $V$ in $\Xi$. On $\Xi$ we have a new
configuration space $\Omega=\{+1,\-1\}^{\Z^{d}}$ with elements
$\xi, \zeta, \omega, \ldots$ and we write ${\cal F} = {\cal F}_{\Xi}$.  Of
course, every
$\sigma \in \Omega'$ gives rise to a unique $\xi=\xi(\sigma) \in \Omega$
via $\xi_i = \sigma_{(i,0)}, i\in \Z^d$ and much of the structure and notation of
the spin
system on $\Xi$ is inherited quite straighforwardly from that on $\Z^{d+1}$.
For example, given $\xi \in \Omega$, we put $\xi^V =
\xi(\sigma^\Lambda) = (\xi(\sigma))^V$ for $\Lambda \cap \Xi =
V$ and $\sigma_{(i,0)}=\xi_i$.  We also write $\zeta_V\xi_{V^c}$ for the
configuration which is equal to $\zeta$ on $V$ and is equal to
$\xi$ on $V^c$ when given $\xi, \zeta \in \Omega$.

This paper is
about the restriction of the Ising states
$\mu_{\beta}^{h}$ and $\mu_\beta^\pm$ to this layer $\Xi$.  In other words,
with $(X(x), x\in \Z^d)$ the random field corresponding to the
considered $d+1$- dimensional Ising state, we want to study the $d$-
dimensional random field $Y$ with
\begin{equation} \label{resran}
Y(i) = X((i,0)), i \in \Z^{d}.
\end{equation}
Obviously, the distribution of $Y$ is the one induced from that of
$X$.
Writing $\nu_{\beta}^{h}$ and $\nu_\beta^\pm$
(or, $\nu_\beta$ in general) for this induced law from, respectively, the
$\mu_{\beta}^{h}$ and $\mu_\beta^\pm$ we have for example
\begin{equation} \label{resmea}
\int f(\xi) d\nu_\beta^+(\xi) \equiv \nu_\beta^+(f) = \mu_\beta^+(f) \equiv
\int f(\sigma) d\mu_\beta^+(\sigma) \end{equation}
for the expectation of any function $f$ which is $\cal F$- measurable
(depends only on the
$\sigma_{(i,0)}, i \in \Z^{d}$).  In particular, for $\beta > \beta_c$ we have
$\nu_\beta^+(Y(0)) = \mu_\beta^+(X(0)) \equiv m^\star(\beta) > 0$ (the
spontaneous magnetization).  Similarly, the truncated correlations (or
covariances) within the layer $\Xi$
\begin{equation} \label{corr}
\nu_\beta(f,g) \equiv
\nu_\beta(fg) -
\nu_\beta(f)
\nu_\beta(g) =
\mu_\beta(f,g) \equiv
\mu_\beta(fg) -
\mu_\beta(f)
\mu_\beta(g)
\end{equation}
between any two functions $f$ and $g$ depending on a finite number of
coordinates in $\xi$, decay exponentially fast in the distance between the
dependence sets of $f$ and $g$, whenever this is the case in the
considered $d+1$- dimensional Ising state $\mu_\beta$ (which is verified
away from
the critical point $\beta=\beta_c$).\newline
The problem can therefore
not be
to evaluate the expectation value of specific observables in our restricted
state because this can be done starting from the well-known Ising states.
Rather, we are interested in some global characterizations of
the restricted states $\nu_\beta$.  More specifically, we wish to
understand the $\nu_\beta$ as Gibbs measures for some interaction.

In the study of the convergence properties of this 
interaction potential it will turn out to be useful to
know whether for the original Ising measure $\mu_\beta$ 
there is good decay of correlations close to the
surface when on the surface we impose a typical configuration drawn from 
$\nu_\beta$. This is the context of the following proposition. 
For $\mu_\beta$ the
Ising measure  on $\{-1,1\}^{\Z^{d+1}}$,
define
\begin{equation}
\label{murand}
\mu_\beta^\xi (\cdot) := \mu_{\beta}(\cdot|\mathcal{F})(\xi).
\end{equation}
This is defined for $\nu_\beta$ almost every surface configuration 
$\xi  \in \Xi$.

One has the following result

\begin{proposition}
\label{prop1}
Suppose that there are constants $C < \infty, m > 0$ so that for all
$x\in \Z^{d+1}, \mu_\beta(X(0) ; X(x)) \leq C \exp(-2m |x|)$.
Then there is a set $K_0\subset \Omega$ with $\nu_\beta(K_0)=1$ such that for all
$\xi \in K_0$ there is a length $\ell=\ell(0,\xi) < \infty$ for which

\begin{equation} 
\label{dec}
|\mu_\beta^\xi(X(0,1) ; X(i,1))| \leq C e^{-m |i|}
\end{equation}
whenever $|i| > \ell$.
\end{proposition}
{\bf Proof} :  The crucial step is to observe
that (by definition) $\nu_\beta( \mu_\beta^\xi) =
\mu_\beta$. Moreover, by the FKG-inequality (positive correlations),
$\nu_\beta[\mu_\beta^\xi(X(0))\mu_\beta^\xi(X(x))]
\geq \nu_\beta[\mu_\beta^\xi(X(0))]\nu_\beta[\mu_\beta^\xi(X(x))] =
\mu_\beta(X(0)) \mu_\beta(X(x))$ since $\mu_\beta^\xi(X(x))$ is
a bounded measurable function non-decreasing in $\xi$.  The constrained
measure $\mu_\beta^\xi$ is itself an FKG-measure so that
$\nu_\beta(|\mu_\beta^\xi(X(0,1) ; X(i,1))|) =
\nu_\beta[\mu_\beta^\xi(X(0,1) ; X(i,1))] \leq \mu_\beta(X(0,1);X(i,1))
\leq
C \exp(-2m |i|)$. The conclusion (\ref{dec}) now follows from standard
Borel-Cantelli 
arguments.
\QED

We feel that the left hand side of (\ref{dec}) is monotone non-increasing
in $\xi$ (when we put more plusses) at least when $h\geq 0$ and $\xi$ has
positive average magnetization.  Proving that $\ell(\xi)$ is decreasing in $\xi \in K_0$
seems to ask however for a rather non-trivial extension of a GHS-type
inequality, \cite{GHS}

\subsection{Results on Gibbsian characterizations}

As announced in Section 1, we restrict ourselves to results concerning
Gibbsian descriptions of restrictions to a layer of Ising states. We first
give a summary of results
describing the state of the art before Dobrushin's 1995 talk, \cite{D}.
We then present the results of the Dobrushin program for these restrictions
of the Ising model.

The beginning of the study of Ising restrictions was

\begin{theorem}[Schonmann, \cite{Sch}]\label{rhs}
Take $\beta>\beta_c$.  The projection $\nu_\beta^+$ of the
two-dimensional plus
phase $\mu_\beta^+ (d=1)$ is non-Gibbsian : there is no translation
invariant absolutely and uniformly convergent potential for
$\nu_\beta^+$.
\end{theorem}

In the following, a further decimation $\nu_\beta^{+,b}, b=3, 4, \ldots$ of
this
$\nu_\beta^+$ was considered.  This measure is obtained as the restriction
of the two-dimensional $\mu_{\beta}^+$ to $\{0\} \times b\Z$ or,
alternatively, as the restriction of $\nu_\beta^+$ to the decimated
integers $b\Z$.

\begin{theorem}[L\H{o}rinczi, Vande Velde, \cite{LVV}]\label{lv}
For sufficiently large $\beta$, for all $b = 3, 4, \ldots,
\nu_\beta^{+,b}$ is a (bona-fide) Gibbs measure.
\end{theorem}

Decimation of non-Gibbsian measures can thus be ({\it bona fide}) Gibbsian
measures (and the opposite is also true). We can extend this result to
random decimations. In other words,
we assign a Bernoulli variable $n_i=0,1$ to each site $i\in \Z$.  The $n_i$
are independent and identically distributed with density $p$. We
consider the restriction $\nu_\beta^{+,(n_i)}$ of $\mu_\beta^+$ (or
$\nu_\beta^+$) to the (random) set $\{i\in \Z : n_i = 1\}$ of occupied
sites.

\begin{theorem} \label{nunu}
There is $p_o > 0$ so that for sufficiently large $\beta$, for all $p <
p_o, \nu_\beta^{+,(n_i)}$ is a (bona-fide) Gibbs measure for almost all
$(n_i)$.
\end{theorem}

One can ask what happens when the temperature is large or when the magnetic
field is non-zero.  While it is rather easy to show that for sufficiently
small $\beta > 0$ or for sufficiently large $h$, the Ising restrictions
$\nu_\beta^h$ are Gibbsian, it is less trivial to show the following

\begin{theorem}[L\H{o}rinczi, \cite{EL}]\label{el}
For $d=1$, the Ising restriction $\nu_\beta^h $ is Gibbsian
whenever $h\neq 0$.
\end{theorem}

Theorem \ref{rhs} was given a more intuitive proof  in \cite{EFS}.  In
fact, something more was obtained (so called absence of quasi-locality).

\begin{theorem}[van Enter, Fernandez, Sokal, \cite{EFS}]\label{efs}
Take $\beta$ sufficiently large.  For any specification $\Gamma$  with
$\nu_\beta^+ \in {\cal G}(\Gamma), d=1$
$\nu_\beta^+$ is not a (bona fide) Gibbs measure.
\end{theorem}

This was extended to any dimension by

\begin{theorem}[Fernandez, Pfister, \cite{FP}]\label{fp}
Take $\beta$ sufficiently large.  For any specification $\Gamma$  with
$\nu_\beta^+ \in {\cal G}(\Gamma), d\geq 1$,
$\nu_\beta^+$ is not a (bona fide) Gibbs measure.
\end{theorem}

On the positive side and as we already mentioned following property 4. in
Definition \ref{specification}, from \cite{FP} it also follows there exists (for
all $\beta,h$)  an everywhere right-continuous local specification $\Gamma =
\Gamma_\beta^h$ 
so that $\nu_\beta^h \in {\cal G}(\Gamma)$.  We observed via
(\ref{vac-gamma}) (see also \cite{MRV} that this implies that the
corresponding vacuum potential
is always convergent (on $\Omega$)).  As we have shown again around
(\ref{vac-gamma}), this implies that $\nu_\beta^h$ (including $h=0^+$) is
weakly Gibbsian for the (everywhere) convergent vacuum potential.
The question about the absolute convergence of the vacuum potential (on
a set of measure one) was also solved:

\begin{theorem}[Dobrushin, Shlosman, \cite{DS}]\label{ds}
For $d=1$ and $\beta$ sufficiently large, $\nu_\beta^+$ is weakly
Gibbsian for the absolutely convergent vacuum potential.
\end{theorem}

For the telescoping potential we have the following

\begin{theorem}[Maes, Vande Velde, \cite{MV}]\label{mv}
For $d=1$ and $\beta$ sufficiently large, $\nu_\beta^+$ is weakly
Gibbsian for the absolutely convergent telescoping potential.
\end{theorem}

The next Section will start with a more detailed presentation (and proof)
of this last Theorem.

\section{Proofs} 
\setcounter{equation}{0}

We use here the telescoping potential 
constructed in Section 3 to prove the results
of Section 4. The main thing to show is the exponential decay
of this potential for large sets which will follow
from the fact that it can be expressed as a correlation function
in a two dimensional Ising model on a halfplane with a "typical"
surface configuration. The specification $\Gamma$ used
in this section will always be the monotone right-continuous specification
consistent with the restriction of the plus phase of the two dimensional
Ising model to the line $\{ (i,0):i\in \Z \}$.
The telescoping potential introduced in Section 3 is here
for $j \not= k$
\begin{equation}
\label{onedimpot} U([j,k],\xi) = - \ln \frac {\gamma_{[j,k]}(\xi^{[j,k]}|+)
\gamma_{[j,k]}(\xi^{]j,k[}|+)}
{\gamma_{[j,k]}(\xi^{]j,k]}|+)
\gamma_{[j,k]}(\xi^{[j,k[}|+)}.
\end{equation}
Using that $\Gamma$ is a specification consistent with the restriction of the
plus phase of the two-dimensional Ising model, we have 
\begin{eqnarray} 
\label{1dimpo}
-U([j,k],\xi) & = & \frac 1{2} (1 - \xi_j)(1 - \xi_k) \ln \frac
{\mu_\beta^{+,\xi^{[j,k]}}[e^{2\beta X(j,1)} e^{2\beta X(k,1)}]}
{\mu_\beta^{+,\xi^{[j,k]}}[e^{2\beta X(j,1)}]
\mu_\beta^{+,\xi^{[j,k]}}[e^{2\beta X(k,1)}]} \nonumber \\
& + & \beta(1-\xi_j)(1-\xi_{j+1}) \delta_{j,k-1},
\end{eqnarray}
where $\mu_\beta^{+,\xi^{[j,k]}}$ is the constrained measure of
(\ref{murand}).  More specifically, consider the events
$S_n(\xi_{[j,k]}) \equiv \{X(x) = +1, x\in
\Lambda_n^c ; X(x) = \xi_i, x=(i,0), i \in [j,k] ; X(x) = +1, x=(i,0), i \notin
[j,k]\}$ where $\Lambda_n$ is an increasing sequence of squares centered
around the origin. For a continuous function $f$
on $\Omega^\prime$,
\begin{equation} \label{mamamia}
\mu_\beta^{+,\xi^{[j,k]}}[f] \equiv \lim_n \mu_\beta^+[f|S_n(\xi_{[j,k]})]
= \lim_n \mu_{n,\beta}^{+,\xi^{[j,k]}}(f).
\end{equation}
This, of course, is a function of the $\xi_i, i\in [j,k]$ only, which act
as extra boundary conditions. The limit (\ref{mamamia}) is over the finite
volume Ising measures $\mu_{n,\beta}^{+,\xi^{[j,k]}}$ with plus boundary
conditions outside the square
$\Lambda_n$ and $\xi^{[j,k]}$ boundary conditions in the middle of the
square (on $\Xi \cap \Lambda_n$, cutting the square in two equal parts).  For
$j=k$, we have $U([j,k],\xi) =$
\begin{equation} \label{mag1}
U(j,\xi) = (1 - \xi_j) \ln \mu_\beta^{+,\xi^j}[e^{2\beta X(j,1)}]
 + \beta(1-\xi_j).
\end{equation}
The extra term for $j=k-1$ and $j=k$ in
(\ref{1dimpo}) and (\ref{mag1}) comes from the interaction 
inside the layer $\Xi$ $(\sim \Z$) and corresponds to the one-dimensional
Ising model. 
One can check that, $\max_\xi |U_A(\xi)| \leq
10\beta $ uniformly in $A$.

Let $\Omega_U
= \cap_{i\in \Z} \Omega_U^i$ be defined via
\begin{eqnarray}
\label{go1}
\Omega^i_U & = & \{\omega\in \{+,-\}^{\Z} : \ \exists
\ell^+_i(\omega),\ell^-_i(\omega) < \infty, \nonumber \\
 & & \forall k > \ell^+_i(\omega),
\forall l > \ell^-_i(\omega), \nonumber \\
& & 1/k\sum_{j=0}^k \omega(i+j) > 8/9,
\quad 1/l\sum_{j=-l}^{-1} \omega(i+j) > 8/9\}.
\end{eqnarray}
(This notation suggests of course that $\Omega_U$ coincides with the 
tail field set 
$\Omega_U$ of points of absolute convergence of the potential $U$
introduced in Section 2. We show that this is indeed the
case.)

Clearly, $\Omega_U \in  \mathcal{T}^{\infty}$ 
is a translation invariant set in the tail field.
It is easy to see that $\nu_\beta^+(\Omega_U)=1$
whenever $\beta$ is sufficiently large.  In fact, under this
condition,
for $\nu_\beta^+$, the  $\ell^{\pm}_i(\cdot), i \in \Z$ are exponential
random variables.
In the following proposition we show that the potential $U$ satisfies
the bound (\ref{hope}) of section 3. For the sake of completeness we repeat
the proof of \cite{MV}

\begin{proposition} \label{cmp}
The potential
defined above (see (\ref{1dimpo}), (\ref{mag1}))
is absolutely convergent for all $\xi \in \Omega_U$
(see (\ref{go1})). In particular, there
are constants $C=C_\beta < \infty$,
$\lambda = \lambda(\beta) > 0$ so that for all $\xi \in \Omega_U, k \in \Z$
\begin{equation} \label{crux}
|U([j,k],\xi)| \leq C e^{-\lambda |j-k|}
\end{equation}
whenever $|j-k| > \ell^-_k(\xi)$. The assumption (\ref{tata}) holds or,  
for all
$\xi \in \Omega_U, i\in \Z$, there is a constant $c_\beta(i,\xi) < \infty$
\begin{equation}
\sum_{k\geq i} \sum_{j\leq i} |U([j,k],\xi)| \leq
16\beta \sum_{k\geq i : \ell^-_k(\xi) > k-i} \ell^-(k,\xi) +
C(\frac{e^\lambda}{e^\lambda - 1})^2 \leq c_\beta(i,\xi).
\end{equation}
\end{proposition}

{\bf Proof:}
Looking at
(\ref{1dimpo}), we see that the crux of the matter consists 
in proving that 
uniformly in the size
$n$ of the boxes $\Lambda_n$
\begin{equation}   \label{must}
\mu_{n,\beta}^{+,\xi^{[0,k]}}[X(0,1);X(k,1)] \leq C e^{-\lambda k}
\end{equation}
whenever $n > k > \ell^+(0,\xi)$.
This was done in
\cite{MV}. 

We repeat the two main steps. They were inspired by the proof of some ergodic
properties of the plus phase in \cite{BS}.
The first step reformulates the required estimate in terms of a percolation
event. Denote by $E_n(o,k)$ the event that there is a path of consecutive
nearest neighbor sites 
$x=(i,j) \in \Lambda_n$, $j>0$, connecting $x=(0,1)$ with $x=(k,1)$ on which
$ (X(x),X'(x)) \not= (+1,+1)$.
Here $X$ and $X'$ are two independent copies of the random field with law
$\mu_{n,\beta}^{+,\xi^{[0,k]}}$.
Then,
\begin{equation}
\left|
\mu_{n,\beta}^{+,\xi^{[0,k]}} \left(X(0,1);X(k,1)\right)
\right| \leq
 2 \mu_{n,\beta}^{+,\xi^{[0,k]}} 
 \times  \mu_{n,\beta}^{+,\xi^{[0,k]}} [E_n(0,k)].
\end{equation}

For the second step, we use  that there is a finite constant $C$
so that 
\begin{equation}
\mu^+_{n,\beta} \times \mu^+_{n,\beta} \left[ E_n(0,k) \right]
\leq C e^{-2\beta k}
\end{equation}
for all sufficiently large $\beta$, uniformly in the size $n$
(see \cite{BS}). 
The argument is now completed by noticing that
\begin{equation}\label{step}
\mu_{n,\beta}^{+,\xi^{[0,k]}} 
 \times  \mu_{n,\beta}^{+,\xi^{[0,k]}}
 [E_n(0,k)]
 \leq
 \exp\left[ 4\beta \sum_{i\in [0,k]} (1-\xi_i) \right]
 \mu_{n,\beta}^{+} 
 \times  \mu_{n,\beta}^{+} [E_n(0,k)].
\end{equation}
The conclusion is that 
\begin{equation}
\label{conclusion}
\left|
\mu_{n,\beta}^{+,\xi^{[0,k]}} \left(X(0,1);X(k,1)\right) \right|
\leq C' e^{4\beta \sum_{i=0}^k(1-\xi_i)} e^{-2\beta k}.
\end{equation}
If $\xi \in \Omega_U$, only one spin out of eight can be minus in $[0,k]$
for $k>l^+(0,\xi)$ and hence the right hand side of (\ref{conclusion})
is then smaller than $e^{4\beta\frac{1}{8}k} e^{-2\beta k}$.
\QED

{\bf Remark 1:} Notice that no use was made of a cluster expansion in the
proof above.  In fact, a naive application of this cluster expansion is
quite impossible as it would yield too much; the attempt in \cite{MV92}
failed for that reason.  This is similar to the
analysis of Gibbs fields for a random interaction in the Griffiths'
regime, see e.g. \cite{BKL}, \cite{vDKP}, \cite{GM}. We must only
concentrate on a specific covariance and percolation techniques seem to be
rather powerful in such cases.

{\bf Remark 2:}
An important ingredient in the previous proof is in the step (\ref{step}).
Therefore it seems that the proof is necessarily restricted to one
dimension.  This however is not the case.  We can prove the decay of the
covariance (as in (\ref{must})) also in higher dimensions. This we will
deal with in a future publication. 

So far we dealt with Theorem 4.8. We now prove the other theorems
of Section 4. For the regular decimation we can only prove  
Theorem \ref{lv} for $b>4$ whereas L\H{o}rinczi {\it et al.} included also 
$b=3,4$.   

\vspace{3mm}

\underline{a)The case $h\neq 0$ or $T>T_c (\beta < \beta_c)$ (Theorem 4.4).}

In this case we don't need the steps above.  The covariance (the left hand
side of (\ref{conclusion})) is exponentially small uniformly in the boundary
condition.  That is,
\begin{equation}
|\mu _n^{+,\xi }( X(0,1) X(k,1) )-\mu _n^{+,\xi }(X(0,1))\mu
_n^{+,\xi} (X(k,1))| \leq C \exp[-\lambda k].
\end{equation}
This is an immediate application of the result that for $h\neq
0$ or for $\beta \leq \beta_c$ the two-dimensional Ising model has a
completely analytic interaction, \cite{SS}.  Therefore the
telescoping potential is actually uniformly absolutely convergent
and thus in this case the projection is a ({\it
bona fide}) Gibbs measure. 

\underline{b) Regular decimation (Theorem 4.2)).}

In this case we must follow the above analysis with
the only change that (\ref{1dimpo}) must be slightly changed. We now have

\begin{eqnarray} 
\label{1dimpodec}
\lefteqn{
-U([j,k],\xi) =  
\frac 1{4} (1 - \xi_j)(1 - \xi_k) } \nonumber \\
&&\ln \frac
{\mu_\beta^{+,b,\xi^{[j,k]}}[e^{2\beta [X(j,1)+X(j,-1)]} e^{2\beta
[X(k,1)+X(k,-1)]}]} {\mu_\beta^{+,\xi^{[j,k]}}[e^{2\beta [X(j,1)+X(j,-1)]}]
\mu_\beta^{+,\xi^{[j,k]}}[e^{2\beta [X(k,1)+X(k,-1)]}]},
\end{eqnarray}
where $\mu_\beta^{+,b,\xi^{[j,k]}}$ is a new constrained measure (very
similar to (\ref{murand})).  More specifically, consider the events
$S^b_n(\xi_{[j,k]}) \equiv \{X(x) = +1, x\in
\Lambda_n^c ; X(x) = \xi_i, x=(i,0), i \in [j,k]\cap b\Z ; X(x) = +1, x=(i,0),
i \in
[j,k]^c \cap b\Z\}$ where $\Lambda_n$ is an increasing sequence of squares
centered around the origin. For a continuous function $f$
on $\Omega^\prime$,
\begin{equation} \label{mamamiadec}
\mu_\beta^{+,b,\xi^{[j,k]}}[f] \equiv \lim_n
\mu_\beta^+[f|S^b_n(\xi_{[j,k]})] = \lim_n
\mu_{n,\beta}^{+,b,\xi^{[j,k]}}(f). \end{equation}
The final trick
of  (\ref{step}) can be repeated but now, since the $\xi$'s live on a
decimated lattice, there are even fewer minusses in the interval $[0,k]$.
It suffices that $|b\Z \cap [0,k]| < k/4$ to have an exponential
decay of the appropriately modified covariance (\ref{must}) uniformly in
the $\xi$, and thus the decimations are Gibbs for $b>4$. 

\underline{c) Random decimations.}

The analysis is as in the previous case.  The final trick of (\ref{step})
must now consider
\begin{equation}
4\beta \sum_{i\in [0,k]} n_i (1-\xi_i).
\end{equation}
Uniformly in $\xi$ this is smaller than $8\beta pk$ for $k$ large
enough on a set of  Bernoulli
variables $(n_i)$ of full measure.  It is
therefore sufficient to choose the density $p < 1/4$.
Note that the variables $n_i$ do not even need to be independent, they
can be chosen according to an ergodic measure $\lambda$ with
$\lambda (n_i ) < 1/4$.
\QED

\section{Variational Principle}
\setcounter{equation}{0}

\subsection{Existence of thermodynamic functions}
\label{var:func}
\subsubsection{Energy density}

We start from the telescoping potential $U$ defined in (\ref{1dimpo}).
In this section  
$\Gamma:=\{ \gamma_V (\sigma|\omega ): V\subset \Z \}$ will as always 
denote the right-continuous specification
such that  the restriction of the plus phase of the two dimensional
Ising model $\nu_\beta^+ \in \mathcal{G} (\Gamma )$. The set 
$\Omega_U \subset \Omega$ will always denote the tailfield set introduced
in (\ref{go1}) on which the telescoping potential is absolutely convergent. 
Given a configuration
$\sigma\in \Omega_U$ we define
\begin{equation}
f_U(\sigma) := \sum_{A \ni 0} \frac{1}{|A|} U_A(\sigma).
\end{equation}
Given a probability measure $\mu$ on $\Omega$ such that 
$\mu (\Omega_U) = 1$ and $f_U\in L^1(\mu )$, we define
\begin{equation}\label{ee}
e^U_\mu := \int f_U(\sigma) d\mu(\sigma).
\end{equation}
Next we introduce the interaction energy in a finite volume $V \subset \Z$:
\begin{itemize}
\item
free boundary conditions:
\begin{equation}
H^f_V(\sigma) := \sum_{A\subset V} U(A,\sigma)
\end{equation}
\item
boundary condition $\omega$:
\begin{equation}
H^\omega_V (\sigma) := \sum_{A \cap V \neq \emptyset} 
U(A,\sigma_{ V} \omega_{ V^c}).
\end{equation}
\end{itemize}
The last sum is well-defined whenever $\omega \in \Omega_U$. For the
potential $U$ constructed in Section
\ref{section3}
we show in this subsection 
that the expectation of the (interaction) energy density exists
and equals $e_\mu^U$ for a certain class of
measures and boundary conditions.

In what follows 
we use $\mathcal{T}$ to denote the set of translation invariant
probability measures on $\Omega$.
The symbols $C,K,c,\lambda$ 
will always be constants whose values can vary from place
to place. We still need the following definitions (cf. Section 5):
\begin{enumerate}
\item
 For $i\in \Z$, $1\leq \alpha < 9/8$, $\sigma \in \Omega$, put
 \begin{equation}
 \label{l+}
 l_i^{\alpha,+} (\sigma) := {\rm min} \{n \in \N : \forall k\geq n
 \quad \frac{1}{k} \sum_{j=0}^{k-1} \sigma_{i+j} \geq \alpha\frac{8}{9} \}
 \end{equation}
 and
 \begin{equation}
 \label{l-}
 l_i^{\alpha,-} (\sigma) := {\rm min} \{n \in \N : \forall k\geq n
 \quad \frac{1}{k} \sum_{j=0}^{k-1} \sigma_{i-j} \geq \alpha\frac{8}{9} \} .
 \end{equation}
 
One must identify 
$l^{1,+/-}_i (\sigma)$ 
with the (abstract) $l_i(\sigma)$ introduced before 
(cf.  Proposition \ref{prop1} and the second remark after Proposition
\ref{prop2}.)
\item 
Let $\mathcal{M}_\alpha$ denote the set of probability measures
$\mu\in\mathcal{T}$ for which $\mu [ l^{\alpha,+/-}_0 (\eta ) >n]\leq e^{-cn}$,
$n\in \N$, for some $c>0$,
and
\begin{equation}
\label{MathcalM}
\mathcal{M} := \bigcup_{\alpha > 1} \mathcal{M}_\alpha.
\end{equation}
Notice that for $\mu \in \mathcal{M}$, $e^U_\mu$ of (\ref{ee}) is
well-defined. 
\item
Finally we put, 
 \begin{eqnarray}
 \lefteqn{
 \Omega_\alpha := \{ \omega \in \Omega \  \exists \epsilon > 0, \ 
 \exists N\in \N :} \nonumber \\ 
 && \ \forall i\in \Z \mbox{ with } |i|\geq N,
 \quad l^{\alpha,+/-}_i (\omega) \leq |i|^{ \frac{1}{3+\epsilon} } \}.
 \end{eqnarray}
\end{enumerate}
Note that $\mu \in \mathcal{M_\alpha} \Rightarrow \mu(\Omega_\alpha) = 1$.
Also if $\alpha < \beta$ then $\Omega_\alpha \supset \Omega_\beta$ and 
$\forall \alpha >1$, $\Omega_\alpha \supset \Omega_U$.
The class of measures $\mathcal{M}$ has to be thought of as the analogue of the
class of tempered measures in the context of unbounded spin systems (cf. 
\cite{LP} Definition 4.1).

\begin{proposition}[Energy density for free boundary conditions] 
\label{var:prop1}
Let $\mu \in \mathcal{M}$ and $U$ the potential defined in (\ref{ja1}).
Then,
\begin{equation}
e^U_\mu = {\rm lim}_{V\uparrow \Z} \frac{1}{|V|} \mu( H^f_V).
\end{equation}
\end{proposition}

\begin{proposition}[Energy density for fixed boundary condition]
\label{var:prop2}
Let $\mu \in \mathcal{M}_\alpha$
and $\omega \in \Omega_\alpha$ for some $\alpha > 1$.
Then, 
\begin{equation}
e^U_\mu = {\rm lim}_{V\uparrow \Z} \frac{1}{|V|} \mu(H^\omega_V).
\end{equation}
\end{proposition}

{\bf Proof of Proposition \ref{var:prop1}:}\\
Consider a sequence of intervals
$V_n = \{-n, \ldots,n \}$.
Then by translation invariance of $\mu$
\begin{eqnarray}
\label{energy1}
\lefteqn{
{\rm lim}_{V\uparrow \Z} \left| \left[
 \frac{1}{|V|} \mu (H^f_V) - e^U_\mu \right] \right| } \nonumber \\
 & = & {\rm lim}_{n \rightarrow \infty} \frac{1}{2n+1} 
  \left| \mu \left[ \sum_{A \subset V_n} U(A,\sigma) - 
   \sum_{i\in V_n} \sum_{A \ni i} \frac{1}{|A|} U(A,\sigma) \right] \right| 
    \nonumber \\
 & \leq &   {\rm lim}_{n \rightarrow \infty} \frac{1}{2n+1} 
  \left| \mu \left[ 
   \sum_{A\cap V_n \neq \emptyset, A \cap V_n^c \neq \emptyset} 
    U(A,\sigma) \right] \right| 
\end{eqnarray}
From Section \ref{section3} we know that 
\begin{itemize}
\item
 $U(A,\sigma) = 0$ if $A$ is not an interval
\item
 we have the upperbounds
 \begin{equation}
 \label{Ubound}
 \left| U([i,j],\sigma) \right| \leq C.I[ l_i^+(\sigma) \geq j-i] + 
  C' {\rm e}^{-\lambda(j-i)}
 \end{equation}
 \begin{equation}
 \label{Ubound2}
 \left| U([i,j],\sigma) \right| \leq C.I[ l_j^-(\sigma) \geq j-i] + 
  C' {\rm e}^{-\lambda(j-i)}
  \end{equation} 

 (cf. Proposition \ref{cmp}).
\end{itemize}

Inserting the exponential term of the RHS of (\ref{Ubound}) 
in the right hand side of 
(\ref{energy1})
gives
\begin{equation}
{\rm lim}_{n \rightarrow \infty} \frac{C'}{2n+1} \left(
\sum_{i=-\infty}^n \ \sum_{j=n+1}^\infty \ + \ \sum_{i=-\infty}^{-n-1}
\ \sum_{j=-n}^n \right) {\rm e}^{-\lambda(j-i)} = 0
\end{equation}
and we are left with two sums, of which we treat only the first one,
the second one can be done in an analogous way. We abbreviate in what
follows $l_i:=l_i^+$ and $l_i^\alpha:=l_i^{\alpha,+}$.
\begin{eqnarray}
\lefteqn{
{\rm lim}_{n \rightarrow \infty} 
\frac{C}{2n+1}  \mu \left(
\sum_{i=-\infty}^n \ \sum_{j=n+1}^\infty 
I[ l_i(\sigma) \geq j-i]
\right)  } \nonumber \\
& = & 
  {\rm lim}_{n \rightarrow \infty} 
\frac{C}{2n+1} 
\sum_{i=-\infty}^n \ \mu \left( 
 [l_i(\sigma) - (n-i)] . I[ l_i(\sigma) \geq n+1-i] \right) \nonumber \\
& = & 
  {\rm lim}_{n \rightarrow \infty} 
\frac{C}{2n+1} 
\sum_{i=-\infty}^n \ \sum_{M=1}^\infty M 
\mu \left( 
 l_i(\sigma) = M+n-i \right) \nonumber \\
& \leq &
   {\rm lim}_{n \rightarrow \infty} 
\frac{C}{2n+1} 
\sum_{i=-\infty}^n \ \sum_{M=1}^\infty M 
 {\rm e}^{-(M+n-i)c} \nonumber \\
& = & 
0.
\end{eqnarray}
\QED

{\bf Proof of Proposition \ref{var:prop2}:}
\begin{eqnarray}
\lefteqn{
 \mbox{lim}_{V\uparrow \Z} \left[
 \frac{1}{|V|} \ \mu (H^\omega_V) - e^U_\mu \right]  \quad = \quad
{\rm lim}_{V\uparrow \Z} \quad
 \frac{1}{|V|} }\times  .  \nonumber \\
 & & \mu 
 \left(  \sum_{A\cap V \neq \emptyset, A\cap V^c \neq \emptyset}
    \left[ U(A,\sigma_{ V} \omega_{V^c})
    - U(A,\sigma) \right] \right). 
\end{eqnarray}
The second term goes to zero as $V\uparrow \Z$ by the proof of Proposition
\ref{var:prop1}.
Again we take intervals $V_n = \{-n, \ldots, n\}$ and we use the bound
(\ref{Ubound}).
The contribution of the exponential part goes to zero so we only have
to show that
\begin{equation}
\label{energy2}
{\rm lim}_{n \rightarrow \infty} \frac{1}{2n+1} \ 
 \mu_n \left( \sum_{i=-\infty}^n
  [l_i(\sigma_{V_n}\omega_{V_n^c}) - (n-i)] . 
  I[l_i(\sigma_{V_n} \omega_{V_n^c}) > n-i]
   \right)=0.
\end{equation}
Abbreviate $\sigma\omega:=\sigma_{V_n}\omega_{V_n^c}$ and
$f(n):= n^{\frac{1}{3+\epsilon}}$. Fix $\omega\in\Omega_\alpha$,
$\mu\in\mathcal{M}_\alpha$ and $n>0$ large enough such that 
$l_n (\omega ) \leq f(n)$.
We distinguish two cases.
\newline
{\bf Case 1: $l^\alpha_i (\sigma ) \leq (n-i) $:} here we first show
that the set $\{ \sigma \in \Omega: l_i (\sigma\omega ) > (n-i),
\quad l_i^\alpha (\sigma ) \leq (n-i) \}$ is not empty for only
a limited number of $i$'s (for $n$ and $\alpha$ fixed). Indeed,
since $l_i^\alpha (\sigma ) \leq (n-i)$ we have
\begin{equation}
\sum_{k=i}^n \sigma (k) \geq \alpha\frac{8}{9} (n-i+1),
\end{equation}
and thus $l_i (\sigma\omega )>(n-i)$ can only happen if there is a
$p> 1 $ such that
\begin{equation}
\label{pee}
\alpha\frac{8}{9} (n-i+1) + \sum_{k=n+1}^{n+p}\omega (k) < \frac{8}{9} 
(n-i+1+p).
\end{equation}
Therefore such a $p$ cannot be too small:
\begin{equation}
\alpha\frac{8}{9} (n-i+1) - p < \frac{8}{9} (n-i+1+p),
\end{equation}
i.e. $p> K (n-i+1)$ where $K= 8(\alpha -1 )/17 >0$.
On the other hand
 since $\omega\in\Omega_\alpha$, $p\geq f(n+1)$ would imply
\begin{equation}
\sum_{k=n+1}^{n+p} \omega (k) > \alpha \frac{8}{9} p
\end{equation}
therefore
(\ref{pee}) can only be satisfied if $p< f(n+1) $. 
Combining the two inequalities we obtained for $p$ we get

\begin{equation}
(n-i+1) < K^{-1} f(n+1).
\end{equation}
and for  $l_i^\alpha ( \sigma ) \leq (n-i) $ and 
$ \omega \in \Omega_\alpha$, 
\begin{equation}
l_i (\sigma\omega ) \leq (n-i+1) + f(n+1).
\end{equation}
We can thus estimate
\begin{eqnarray}\label{afsch}
\lefteqn{\sum_{i=-\infty}^n \mu_n \left [ \left [ l_i (\sigma\omega )
-(n-i) \right ] I \left [ l_i(\sigma\omega ) > (n-i) \right ] 
| l_i^\alpha (\sigma ) \leq (n-i) \right ] }\nonumber \\
&\leq &
\sum_{i= n - K^{-1} f(n+1)+1}^n \mu_n \left [ l_i (\sigma\omega ) - (n-
i)| l^\alpha_i (\sigma ) \leq (n-i) \right ] \nonumber \\
&\leq &
\sum_{i= n - K^{-1}f(n+1)+1}^n \sum_{j=1}^{f(n)} j 
\quad = \quad O ((f(n))^3).
\end{eqnarray}
{\bf Case 2:} 
$l_i^\alpha (\sigma)> (n-i+1)$: If $\omega\in \Omega_\alpha$,
then, for $p> f(n+1)$,
\begin{equation}
\sum_{j=n+1}^{n+p} \omega (j) > \alpha \frac{8}{9} p.
\end{equation}
Therefore,
\begin{equation}
\sum_{j=i}^n \sigma (j) + \sum_{j=n+1 }^{n+p} \omega (j) \geq -(n-i+1) +
\alpha\frac{8}{9}p.
\end{equation}
Hence 
\begin{equation}
-(n-i+1)+\alpha \frac{8}{9} \geq ((n-i+1) +p) \frac{8}{9}
\end{equation}
implies $l_i( \sigma\omega )\leq p+(n-i+1)$. I.e. if 
$p\geq K^{-1}(n-i+1)$, then $l_i (\sigma\omega ) \leq p+n-i+1$,
and thus we conclude
\begin{equation}\label{ellie}
l_i (\sigma\omega ) \leq (n-i+1)(1+K^{-1}).
\end{equation}
On the other hand, since $\mu\in\mathcal{M}_\alpha$, 
\begin{equation}\label{muu}
\mu (l_i^\alpha (\sigma ) > (n-i)) \leq e^{-c (n-i)}.
\end{equation}
Combining (\ref{ellie}) and (\ref{muu}) we get
\begin{eqnarray}
\label{afsch2}
\lefteqn{
\sum_{i=-\infty}^n \mu_n
\left( ( l_i (\sigma\omega ) - (n-i) )
I \left [ l_i (\sigma \omega ) > (n-i) \right ]
  | l_i^\alpha ( \sigma ) > (n-i)  \right). } \nonumber \\ 
&&
\mu_n \left(l_i^\alpha (\sigma ) > (n-i) \right) 
\leq  \sum_{i=-\infty}^n K (n-i) e^{-c (n-i)}=O(1).
\end{eqnarray}
Conditioning respectively on 
$l_i^\alpha (\sigma ) \leq (n-i) $, $l_i^\alpha (\sigma ) > (n-i)$, 
and using (\ref{afsch}) (\ref{muu}),
(\ref{afsch2}) we arrive at
\begin{equation}
\sum_{i=-\infty}^n \mu_n \left [ l_i(\sigma\omega ) - (n-i)
I \left [ l_i (\sigma\omega ) > (n-i) \right] \right] \leq 
O((f(n))^3),
\end{equation}
and hence (\ref{energy2}) follows.\QED
\subsubsection{Pressure}
\begin{proposition}
\label{var:prop3}
Let $U$ be the potential defined in (\ref{ja1}). The pressure 
(or free energy density)
\begin{equation}
P(U) := {\rm lim}_{V\uparrow \Z} \frac{1}{|V|} {\rm log} Z^f_V
\end{equation}
exists,
where
\begin{equation}
Z^f_V \equiv \sum_{\sigma \in \Omega_V} \exp \left[ -\sum_{A \subset V}
U(A,\sigma) \right]
\end{equation}
is the finite volume ($V$) partition function with free boundary conditions.
\end{proposition}

{\bf Proof:}
\begin{eqnarray}
\sum_{A \subset V} U(A,\sigma) & = & H^f_V (\sigma) \nonumber \\
& = & {\rm log} \frac{\gamma_V(+|+)}{\gamma_V(\sigma^V|+)}.
\end{eqnarray}
Therefore
\begin{equation}
Z^f_V = \frac{1}{\gamma_V(+|+)}.
\end{equation}
The
existence of ${\rm lim}_V \frac{1}{|V|} {\rm log} \gamma_V(+|+)$
follows from standard arguments for the two dimensional Ising model.
Indeed,
\begin{equation}\label{hahaha}
\lim_{n\uparrow\infty}\frac{1}{|V_n|} \log \gamma_{V_n} (+|+)
=\lim_{n\uparrow\infty}\frac{1}{|V_n|}\log
\frac{Z_{\Lambda_n}}{Z_{\Lambda_n\cup\tilde{V}_n}}.
\end{equation}
Here $Z_\Lambda$ denotes the partition function of the two dimensional
Ising model for volume $\Lambda$ and $+$-boundary conditions
outside $\Lambda$, 
$\Lambda_n=\Lambda_{1,n} \bigcup \Lambda_{2,n}$, where
\begin{equation}
\Lambda_{1,n} := \{ (i,j)\in\Z^2: |(i,j)|\le n, \quad j>0 \},
\end{equation}
\begin{equation}
\Lambda_{2,n} := \{ (i,j)\in\Z^2: |(i,j)|\le n, \quad j<0 \},
\end{equation}
and
\begin{equation}
\tilde{V}_n := \{ (i,0): i\in V_n \}.
\end{equation}
The existence of the limit in (\ref{hahaha}) is thus standard and can
be calculated from the cluster expansion. It is equal to minus 
the free energy
density of the two dimensional Ising model plus a surface contribution.
\QED
\begin{proposition}
\label{var:prop4}
Let $P(U)$ be as in Proposition \ref{var:prop3} and define
\begin{equation}
Z^\omega_V := \sum_{\sigma_V} \exp \left[
 - \sum_{A \cap V \neq \emptyset} U(A,\sigma_{ V} \omega_{
 V^c}) \right].
\end{equation}
Then $\forall \alpha > 1$ $\forall \omega \in \Omega_\alpha$
\begin{equation}
P(U) = {\rm lim}_{V\uparrow \Z} \frac{1}{|V|}{\rm log} Z^\omega_V.
\end{equation}
\end{proposition}

{\bf Proof:} \\
\begin{eqnarray}\label{buu}
\lefteqn{
\frac{1}{|V|} \quad {\rm log} \quad \frac{Z^f_V}{Z^\omega_V} =}
 \nonumber \\
&  &  \frac{1}{|V|} {\rm log} \ \mu_V^{\omega,U} \left( \exp \left[ 
 \sum_{A \cap V \neq \emptyset, A \cap V^c \neq \emptyset} 
 U(A,\sigma_{ V} \omega_{ V^c}) \right] \right)
 \leq 0
\end{eqnarray}
where $\mu_V^{\omega,U}$ is the measure introduced in Section 2,
(2.15).
The inequality in (\ref{buu}) follows from the fact $U(A,.) \leq 0 $ 
$ \forall A, \ |A|\geq 2$. This follows from
the expression (\ref{1dimpo}) for the potential, and the positivity of
correlations for monotonic functions.
On the other hand by Jensen's inequality  
(take $V_n = \{-n, \ldots, n\}$)
\begin{eqnarray}
\lefteqn{
\mbox{lim} \ \mbox{inf}_n \ \frac{1}{2n+1} \ {\rm log} \ 
\mu_{V_n}^{\omega,U} \left( \exp \left[
\sum_{A \cap V \neq \emptyset, A \cap V^c \neq \emptyset} 
 U(A,\sigma_V \omega_{V^c}) \right] \right)
 } \nonumber \\
& \geq & 
\mbox{lim} \ \mbox{inf}_n \ \frac{1}{2n+1}  \ 
\mu_{V_n}^{\omega,U} \left( 
\sum_{A \cap V \neq \emptyset, A \cap V^c \neq \emptyset} 
 U(A,\sigma_{ V} \omega_{V^c})  \right) \nonumber \\
& = & 0
\end{eqnarray}
as can be obtained from the proof of Proposition \ref{var:prop2}.
\QED

\subsection{First part of the variational principle}
For every finite volume $V$, every $\omega \in \Omega_U$, every $\mu \in
\mathcal{T}$, the following holds
\begin{eqnarray}
S_V(\mu|\gamma_V(.|\omega)) & = & -S_V(\mu) + \mu(H^\omega_V) + {\rm log}
Z^\omega_V   
\nonumber \\
& \geq & 0
\end{eqnarray}
where 
$S_V(\mu)$ is the entropy of the measure $\mu$, defined as

\begin{equation}
\label{entropy}
S_V(\mu):= 
 - \sum_{\sigma \in \Omega_V} \mu_V(\sigma) {\rm log} 
       \mu_V(\sigma)
\end{equation}
and
$S_V(\mu|\nu)$ is the relative entropy of the measure $\mu$ with respect
to the measure $\nu$, defined as 
\begin{equation}
\label{relentropy}
S_V(\mu|\nu) :=  \sum_{\sigma \in \Omega_V} \mu_V(\sigma) {\rm log} 
\frac{\mu_V(\sigma)}{\nu_V(\sigma)}
\end{equation}
if $\mu$ is absolutely continuous with respect to $\nu$ and
$S_V (\mu|\nu ) = +\infty$ otherwise (we make the convention $0\log 0=0$).
We still need the following notation:
\begin{equation}
\nu \circ \gamma_V (f) := \int d\nu(\omega) \sum_{\sigma_V} \gamma_V
(\sigma_V|\omega) \ f(\sigma_V\omega_{V^c}).
\end{equation}

\begin{theorem}
\label{vac:theo}
\begin{enumerate}
\item
Let $P(U)$ be as in Proposition \ref{var:prop3}. Then
\begin{equation}
\label{var:sup}
P(U) = {\rm sup}_{\mu \in \mathcal{M}} \left[ 
s(\mu) - e^U_\mu \right]
\end{equation}
where $s(\mu):= {\rm lim}_V \ \frac{1}{|V|} S_V (\mu )$
is the entropy density.
\item
Let $\omega \in \Omega_\alpha$ for some $\alpha > 1$ and $\mu \in \mathcal{M}$.
Then
 \begin{itemize}
 \item
 the relative entropy density $s(\mu|U)$ exists
 \begin{eqnarray}
 \label{var:relentr}
 s(\mu|U) & \equiv & 
 {\rm lim}_{V\uparrow \Z} \frac{1}{|V|} S_V(\mu_V|\gamma_V(\cdot|+) \nonumber \\
  & = & {\rm lim}_{V\uparrow \Z} \frac{1}{|V|} 
  S_V (\mu|\gamma_V(\cdot|\omega))
 \end{eqnarray}
 \item
 $\forall \nu \in \mathcal{M}$
 \begin{equation}
 \label{var:relentr2}
 s(\mu|U) = 
 {\rm lim}_{V\uparrow \Z} \frac{1}{|V|} S_V(\mu_V|\nu \circ \gamma_V) .
 \end{equation}
 \end{itemize}  
\item
The supremum in (\ref{var:sup}) is reached for $\mu = \nu^+_\beta$. 
\end{enumerate}

\end{theorem}

{\bf Proof:}\\
$1$. is a consequence of 
$ S_V(\mu|\nu) \geq 0$, $\forall \mu, \nu \in \mathcal{T}$ and
subsection \ref{var:func}.\\
$3$. is a special case of 2 because $\nu^+_\beta \in \mathcal{M}$,
$\nu^+_\beta\circ\gamma_{V_n} = \nu^+_\beta$ and 
$ s(\nu^+_\beta|\nu^+_\beta) = 0$.\\
$2$. The first statement, (\ref{var:relentr}) follows immediately from the previous
subsection. To prove (\ref{var:relentr2}) we need to prove that
\begin{itemize}
\item
 \begin{equation}
 {\rm lim}_n \frac{1}{2n+1} \int d\nu(\omega) \ \mu (H^\omega_n - H^f_n) = 0
 \end{equation}
\item
 \begin{equation}
 {\rm lim}_n \frac{1}{2n+1} \int d\nu(\omega) \ {\rm log} 
 \frac{Z^f_n}{Z_n^\omega} =0.
 \end{equation}
\end{itemize}

To prove that the first limit is zero, we proceed as in the proof of
Proposition \ref{var:prop2}. Let $c_\mu$ and $c_\nu$ be the
constants appearing in the definition of $\mathcal{M}$ for
$\mu$, resp. $\nu$. We estimate
\begin{eqnarray}
\lefteqn{
\frac{1}{2n+1} \sum_{j=-\infty}^n 
\E_{\mu \times \nu} \left( [l_j(\sigma,\omega) -
(n-j+1)] . I[l_j(\sigma,\omega) > n-j+1] \right) } \nonumber \\
& = & \frac{1}{2n+1} \sum_j \sum_{M=1}^\infty M \mu \times \nu 
 \left[ l_j(\sigma,\omega) = M+n-j+1  \right] \nonumber \\
& \leq & 
  \frac{1}{2n+1} \sum_j \sum_{M=1}^\infty M  
  \mu \times \nu 
 \left[ l_j(\sigma,\omega) = M+n-j+1  | l_j^\alpha(\sigma) \leq n-j+1 \right] 
  \nonumber \\
&  + & \frac{1}{2n+1} \sum_j \sum_{M=1}^\infty M
\mu \times \nu 
  \left[ l_j(\sigma,\omega) = M+n-j+1  | l_j^\alpha(\sigma) > n-j+1 \right] .
  {\rm e}^{-(n-j)c_\mu} 
  \nonumber \\
& \leq & \frac{1}{2n+1} \sum_j \sum_{M=K^{-1}(n-j)}^\infty \ M 
 {\rm e}^{-Mc_\nu} \quad +
    \frac{1}{2n+1} \sum_j {\rm e}^{-(n-j)c_\mu} 
    \sum_{M=1}^{K(n-j)} M  \nonumber \\
&&   + \ \frac{1}{2n+1} \sum_j {\rm e}^{-(n-j)c_\mu}
    \sum_{M=K(n-j)+1}^\infty M .
    \nu \left[ l_n^\alpha(\omega) > M \right] 
    \nonumber \\
& \leq & \frac{1}{2n+1} \sum_j {\rm e}^{-K^{-1}(n-j)c_\nu} 
\ + \ \frac{1}{2n+1}
 \sum_j O((n-j)^2) {\rm e}^{-(n-j)c_\mu} \nonumber \\ 
 & & + \ \frac{1}{2n+1}  \sum_j {\rm e}^{-(n-j)c_\mu}
  \sum_{M=K(n-j)+1}^\infty M {\rm e}^{-Mc_\nu} \nonumber \\
& \rightarrow & 0 \quad {\rm when} \ n \rightarrow +\infty.      
\end{eqnarray}
This implies (by the argument we used to prove the existence of $P(U)$) that
\begin{equation}
{\rm lim}_{n \rightarrow \infty} \frac{1}{2n+1} \int d\nu(\omega) \ 
{\rm log} \frac{Z^f_n}{Z^\omega_n} = 0.
\end{equation}
\QED

\subsection{Second part of the variational principle}
The second part of the variational principle characterizes the maximizers of
(\ref{var:sup}) as the  measures consistent with $\Gamma$. Note that any
maximizer $\mu$ of (\ref{var:sup}) satisfies $s(\mu|\nu^+_\beta)=0$.
To conclude $\mu \in \mathcal{G}(\Gamma)$ from $s(\mu|\nu^+_\beta)$
we need an extra technical condition:

\begin{theorem}
\label{var:theo1}
Suppose that $\mu \in \mathcal{M}$ such
that
\begin{equation}
s(\mu| \nu^+_\beta) :=
\lim_{n \rightarrow \infty} 
\frac{1}{2n+1} S_{\Lambda_n}(\mu|\nu_\beta^+) = 0
\end{equation}
and
\begin{equation}
\label{cond}
\lim_{n\rightarrow \infty}
\mu \left( \exp [2\beta \sum_{i = -n}^n (1-\eta(i))] \right) 
2^n \exp[-\beta n] \ = 0
\end{equation}
then $\mu \in \mathcal{G}(\Gamma)$. 
\end{theorem}
We will show that $\nu^+_\beta$ satisfies the hypotheses of Theorem 
\ref{var:theo1}.

 {\bf Proof of the Theorem:}\\
The first part of the proof 
follows
\cite{Geo} p323 
(the variational principle in the regular Gibbs case). 
We have to show that for all $\Lambda \in \mathcal{L}$ and for every local
function $g$
\begin{equation}
\label{consistency}
\mu \gamma_\Lambda (g) \ = \ \mu(g).
\end{equation} 
We show that this equality holds for $\Lambda =\{0\}$ and every local function
$g$; (\ref{consistency}) then follows from translation invariance and the
positivity of $\Gamma$.

Note that for every $ \Delta \in \mathcal{L}$, $S_\Delta(\mu|\nu^+_\beta) < \infty$.
This implies that for every $ \Delta \in \mathcal{L}$ there exists a $
\mathcal{F}_\Delta$-measurable function $f_\Delta \geq 0$ such that
$\mu = f_\Delta .\nu$ on $\mathcal{F}_\Delta$ 
($f_\Delta = \frac{d\mu}{d\nu^+_\beta}|_{\mathcal{F}_\Delta}$).
In \cite{Geo} it is shown that for every $\epsilon >0$ and for every interval 
$I \ni 0$ there exists a set $ \Delta \in \mathcal{L}$ with $I\subset \Delta$ such that 
\begin{equation}
\nu^+_\beta \left( |f_\Delta - f_{\Delta \setminus \{0\}} | \right) \leq 
\epsilon.
\end{equation}

Fix a local function $g$ and let $I$ be an interval such that $I \ni 0$ and $g
\in \mathcal{F}_I$.
Given $\epsilon >0$ and $I$, fix $\Delta$ as above and define
\begin{equation}
\tilde{g}(\omega) := \sum_{\sigma_0 = +/-} g(\sigma_0\omega_{0^c}) 
\gamma_{\{0\}} 
\left(\sigma_0 | \omega_{\Delta \setminus \{0\}} +_{\Delta^c}\right),
\end{equation}
$\tilde{g} \in \mathcal{F}_{\Delta\setminus\{0\}}$.
Then
\begin{eqnarray}
\label{6eps}
\lefteqn{ \left| \mu \gamma_{ \{ 0\} }(g) - \mu (g) \right| } 
 \nonumber \\ 
 &\leq& \mu \left( | \gamma_{\{0\}}(g) - \tilde{g} | \right)
  +
  \left|\mu(\tilde{g}) - \nu^+_\beta (f_{\Delta \setminus \{0\}}
   \tilde{g}) \right| \nonumber \\
 &+&
  \nu^+_\beta \left(f_{\Delta \setminus \{0\}} .
   | \tilde{g} - \gamma_{ \{0\} } g|  \right)
  +
  \left| \nu^+_\beta \left(f_{\Delta \setminus \{0\}} .
   ( \gamma_{ \{0\} } g - g)  \right) \right| \nonumber \\
 &+&
  \|g\|_\infty \   \nu^+_\beta \left( \left| f_{\Delta \setminus \{0\}} -
  f_\Delta \right| \right) 
  +
   \left| \nu^+_\beta(f_\Delta g ) -  \mu(g) \right| . 
\end{eqnarray} 
Since $\tilde{g} \in \mathcal{F}_{\Delta \setminus\{0\}}$ and $g \in
\mathcal{F}_\Delta$ the second and the last term on the right are zero. The
fourth term vanishes because $\nu^+_\beta \in \mathcal{G}(\Gamma)$ and
$f_{\Delta \setminus \{0\}} \in \mathcal{F}_{\{0\}^c}$. The fifth term is
smaller than $\|g\|_\infty. \epsilon$ because of the choice of $\Delta$.
We are left with the first and the third term. In the quasilocal case they do
not cause any trouble because there
$\tilde{g} \rightarrow \gamma_{\{0\}} g$ in sup-norm as $\Delta \rightarrow \Z$, 
i.e. for $\Delta$ large enough, 
$\| \tilde{g} - \gamma_{\{0\}}g \|_\infty < \epsilon$.
This is not the case here. We put $\Delta:= \{-n, \ldots, n\}$ and write
\begin{eqnarray}
\label{}
\lefteqn{
\left| \left(\gamma_{\{0\}}g - \tilde{g}\right) (\omega) \right| = } 
\nonumber \\
&& 
 \left| \sum_{\sigma_0 = +/-} \left[ \gamma_{\{0\}}(\sigma_0|\omega_{0^c}) - 
  \gamma_{\{0\}}(\sigma_0|\omega_{\Delta\setminus 0} +_{\Delta^c}) \right] 
  g(\sigma_0\omega_{0^c})  \right| \nonumber \\
&&
\leq \|g\|_\infty \ \sum_{\sigma_0 = +/-} \left| \gamma_{\{0\}} (\sigma_0|\omega_{0^c})
- \gamma_{\{0\}} (\sigma_0|\omega_{\Delta\setminus0} +_{\Delta^c}) \right|.  
\end{eqnarray}
We have that 
\begin{equation}
 \gamma_{\{0\}}(\sigma_0 = +|\omega_{0^c}) = \ 
 \frac{1}{1 + \exp \left\{ \sum_{A \ni 0} \left[ U(A,-_0 \omega_{0^c}) - U(A,
 +_0 \omega_{0^c}) \right] \right\} }
\end{equation}

\begin{eqnarray}
\lefteqn{
 \gamma_{ \{0\} }\left(\sigma_0 = +|\omega_{\Delta} +_{\Delta^c} \right)  }
 \nonumber  \\ 
 & = &\frac{1}{1 + \exp \left\{ \sum_{A \ni 0} \left[ 
 U(A,-_0 \omega_{\Delta \setminus 0} +_{\Delta^c}) -
 U(A,+_0 \omega_{\Delta \setminus 0} +_{\Delta^c}) \right] \right\} }
\end{eqnarray}
and an analogous expression for $\sigma_0 = -$. 

Use the inequality
\begin{equation}
\left| \frac{1}{1+ \mbox{e}^x} - \frac{1}{1+ \mbox{e}^y} \right|
\leq | x-y|
\end{equation}
to obtain for the first term in (\ref{6eps})
\begin{eqnarray}
\lefteqn{
\mu \left( | \gamma_{ \{0\} } g - \tilde{g}| \right) } \nonumber \\
& \leq &
2 \|g\|_\infty  
\int d\mu(\omega) \left[ 
\left| \sum_{A \ni 0} \left[
U(A,+_0 \omega_{0^c}) - U(A, +_0 \omega_{\Delta \setminus0} +_{\Delta^c}) 
\right] \right|  \right. \nonumber \\
&&
\left.
+ \left| 
\sum_{A \ni 0} \left[
U(A,-_0 \omega_{0^c}) - U(A, -_0 \omega_{\Delta \setminus0} +_{\Delta^c}) 
\right] \right|
\right]
\nonumber \\
& = & 
2 \|g\|_\infty \ 
\int d\mu(\omega) 
\left[ 
\left| \sum_{A \ni 0, A \cap \Delta^c \neq \emptyset}
U(A,+_0 \omega_{0^c}) \right| 
\right. \nonumber \\
&& +
\left.
\left| \sum_{A \ni 0, A \cap \Delta^c \neq \emptyset}
U(A,-_0 \omega_{0^c}) \right| 
\right]
\end{eqnarray}
where in the second line we used that $U(A,\omega_{\Delta} +_{\Delta^c}) =0$
whenever $A \cap \Delta^c \neq \emptyset$.
This expression goes to zero when $n$ tends to infinity for every $\mu \in 
\mathcal{M} $ (cf. the proof of Theorem \ref{vac:theo}).

For the third term in (\ref{6eps}) it suffices to show that
\begin{eqnarray}
\label{third}
\lefteqn{
\limsup_n \int d\nu^+_\beta (\omega) \frac{\mu(\omega_{\Delta \setminus 0})}
 {\nu^+_\beta (\omega_{\Delta \setminus 0})} }
 \nonumber \\
&& 
\left[ 
\left|
\sum_{A \ni 0, A \cap \Delta^c \neq \emptyset} 
U(A, +_0 \omega_{0^c}) \right| \right.
\nonumber \\
&&
  \left. +  \ \left|  
 \sum_{A \ni 0, A \cap \Delta^c \neq \emptyset} U(A, -_0 
 \omega_{0^c}) \right| 
 \right] 
 \ = \ 0.
\end{eqnarray}  
Using that $U(A,.) \equiv 0$ if $A$ is not an interval together
with the bound (\ref{hop1})
\begin{equation}
|U([j,k],\omega)| \ \leq \ C_1 I \left[ l^+_j(\omega) > k-j \right] +
 C_2 I \left[ l^+_j(\omega) \leq k-j \right] \ 
 \mbox{e}^{-\lambda (k-j)}
\end{equation}
we obtain
\begin{eqnarray}
\lefteqn{
\sum_{A \ni 0, A \cap \Delta^c \neq \emptyset} 
\left| U(A, \sigma_0 \omega_{0^c})  \right| 
 \leq 
  \left( \sum_{j=-\infty}^{-n-1} \ \sum_{k=0}^{+\infty} \ 
+ \ \sum_{j=-n}^0 \ \sum_{k=n+1}^{+\infty} \right) } \nonumber \\
&&
 C_1 I \left[l^+_j(\sigma_0 \omega_{0^c}) > k-j \right] \ + \ 
C_2 I \left[ l_j^+ (\sigma_0 \omega_{0^c}) \leq k-j \right] 
\mbox{e}^{-(k-j)\lambda} .
\end{eqnarray}
For the exponential part we are back in the quasilocal case,
namely the sum goes to zero as $n$ tends to infinity, uniformly in $\omega$.
We continue with $\sigma_0 = +$ and the first sum, i.e. $j < -n$, $k \geq 0$.
The other sum and $\sigma_0 = -$ can be treated in exactly the same way.
\begin{eqnarray}
\lefteqn{
\sum_{j=-\infty}^{-n-1} \ \sum_{k=0}^{+\infty} I \left[ l_j^+(+_0 \omega_{0^c})
> k-j \right] }\\
& = & \sum_{j=-\infty}^{-n-1} \left( l^+_j(+_0 \omega_{0^c}) - |j| \right)
I \left[ l_j^+(+_0 \omega_{0^c}) > |j| \right] \nonumber \\
& =: & \tilde{u}_n(\omega).
\end{eqnarray}
Note that $\tilde{u}_n(.)$ is  monotone non-increasing (cf. (\ref{l+})). 

Denote by $\rho_\beta$ the Bernoulli measure on $\Omega$ with
$\rho_\beta(\sigma_0 = -) = \mbox{e}^{-8\beta}$. Then $\nu^+_\beta \leq \rho$
(see \cite{Holley}). Returning to (\ref{third}) we can write,
using Cauchy-Schwartz
\begin{eqnarray}
\lefteqn{
 \int d\nu^+_\beta(\sigma) \frac{\mu(\sigma_{\Delta \setminus 0})}
  {\nu^+_\beta(\sigma_{\Delta\setminus 0})}
   \tilde{u}_n(\sigma) }\\
 & = &
 \sum_{\sigma_{\Delta\setminus 0}} \nu^+_\beta (\sigma_{\Delta\setminus 0})
   \frac{\mu(\sigma_{\Delta \setminus 0})}
  {\nu^+_\beta(\sigma_{\Delta\setminus 0})}
 \frac{\rho(\sigma_{\Delta \setminus 0})}
  {\rho(\sigma_{\Delta\setminus 0})}
  \ 
  \E_{\nu^+_\beta} (\tilde{u}_n | \mathcal{F}_{\Delta \setminus 0})
   (\sigma_{\Delta \setminus 0}) \nonumber \\
 & \leq & 
 \ 
  \left[
    \sum_{\sigma_{\Delta\setminus 0}} \rho (\sigma_{\Delta\setminus 0})
   \left(
   \frac{\mu(\sigma_{\Delta \setminus 0})}
  {\rho(\sigma_{\Delta\setminus 0})} \right)^2
  \right]^{1/2} .
  \nonumber \\
&  &
  \left[
  \sum_{\sigma_{\Delta\setminus 0}} \rho (\sigma_{\Delta\setminus 0})
   \left(
   \E_{\nu^+_\beta} ( \tilde{u}_n | \mathcal{F}_{\Delta \setminus 0} )
    (\sigma_{\Delta \setminus 0}) 
   \right)^2
  \right]^{1/2}.
\end{eqnarray}
Now
\begin{eqnarray}
\lefteqn{
 \rho(\sigma_{\Delta\setminus 0}) = 
 \prod_{i=-n, i \neq 0}^n 
 \left( 1 - \mbox{e}^{-8\beta} \right)^{\frac{1+\sigma(i)}{2}} \ 
 \left(\mbox{e}^{-8\beta} \right)^{\frac{1 - \sigma(i)}{2}} 
 } \\
 & \geq &
  \left( 1 - \mbox{e}^{-8\beta} \right)^{2n} \ 
  \mbox{e}^{-4 \beta \sum_{i=-n, i \neq 0}^n (1 - \sigma(i))}
\end{eqnarray}
or
\begin{equation}
 \frac{1}{\rho(\sigma_{\Delta \setminus 0})}
 \leq  
 4^n \ 
  \mbox{e}^{ 4 \beta \sum_{i=-n, i \neq 0}^n (1 - \sigma(i))} 
\end{equation} 
for $\beta $ large.
Therefore
\begin{eqnarray}
\left[
\sum_{\sigma_{\Delta\setminus 0}} 
   \frac{\mu(\sigma_{\Delta \setminus 0})^2}
  {\rho(\sigma_{\Delta\setminus 0})}
 \right]^{1/2}
 & \leq &
 \left(
 \sum_{\sigma_{\Delta\setminus 0}}
 \left[ 2^n \ 
 \mbox{e}^{2 \beta \sum_{i=-n, i\neq 0}^n (1 - \sigma(i))} \ 
  \mu(\sigma_{\Delta \setminus 0})
  \right]^2
  \right)^{1/2}
  \nonumber \\
  & \leq & 
  2^n \ 
  \mu \left( \mbox{e}^{2 \beta \sum_{i=-n}^n (1 - \sigma(i))} \right).
\end{eqnarray}  
On the other hand
$ \rho \geq \nu^+_\beta$
implies that $\forall \Lambda \in \mathcal{L}$
\begin{eqnarray}
 \sum_{\sigma_\Lambda} \rho (\sigma_\Lambda)
   \E_{\nu^+_\beta} ( f | \mathcal{F}_\Lambda )
    (\sigma_\Lambda) 
&  \geq  & 
 \sum_{\sigma_\Lambda} \nu^+_\beta (\sigma_\Lambda)
   \E_{\nu^+_\beta} ( f | \mathcal{F}_\Lambda )
    (\sigma_\Lambda)  \nonumber \\
 & = &
    \E_{\nu^+_\beta} ( f)
\end{eqnarray}
for every non-decreasing monotone function $f$.
Or, since $-\tilde{u}_n^2$ is monotone non-decreasing
\begin{eqnarray}
\lefteqn{
\sum_{\sigma_{\Delta \setminus 0}} 
   \rho (\sigma_{\Delta \setminus 0})
   \left[
   \E_{\nu^+_\beta} ( \tilde{u}_n | \mathcal{F}_{\Delta\setminus 0} )
    (\sigma_{\Delta \setminus 0}) 
   \right]^2 } \nonumber \\
&  \leq &
  \sum_{\sigma_{\Delta \setminus 0}} 
   \rho (\sigma_{\Delta \setminus 0})
   \left[
   \E_{\nu^+_\beta} ( \tilde{u}_n^2 | \mathcal{F}_{\Delta\setminus 0} )
    (\sigma_{\Delta \setminus 0}) 
   \right] 
   \nonumber \\
& \leq &
   \sum_{\sigma_{\Delta \setminus 0}} 
   \nu^+_\beta (\sigma_{\Delta \setminus 0})
   \left[
   \E_{\nu^+_\beta} ( \tilde{u}_n^2 | \mathcal{F}_{\Delta\setminus 0} )
    (\sigma_{\Delta \setminus 0}) 
   \right] 
   \nonumber \\
& = &   
  \int d\nu^+_\beta (\sigma) \tilde{u}_n^2 (\sigma).  
\end{eqnarray}     
Putting the pieces together we obtain for the third term in
(\ref{6eps}) 
\begin{equation}
\label{final}
\nu^+_\beta \left(
f_{\Delta \setminus 0} . 
\left| \tilde{g} - \gamma_\Lambda g \right| \right)
\ \leq \ 
 C \ 2^n \ \mu \left(
 \mbox{e}^{2 \beta \sum_{i=-n}^n (1 - \sigma(i))} \right) \ 
 \left[ \nu^+_\beta ( \tilde{u}_n^2 ) \right]^{1/2}.
\end{equation}

\begin{lemma}
\label{lem1}
\begin{equation}
\left[ \nu^+_\beta ( \tilde{u}_n^2 ) \right]^{1/2}
 \ \leq \ 
C \ \mbox{e}^{- \beta n}.
\end{equation}
\end{lemma} 
The proof of Lemma \ref{lem1} is straightforward and uses only that
$\nu^+_\beta \left[ l_j^+ (\omega) > n \right] \leq C \mbox{e}^{-\beta n}$.

Lemma \ref{lem1} together with 
(\ref{final}) prove the Theorem.
\QED

The following Lemma states that for $\mu = \nu^+_\beta$ and $\beta$ large enough
the conditions of
the Theorem are fulfilled (as it should be) and a fortiori that the class of
$\mu$'s for which (\ref{cond}) holds is not empty.  
 
\begin{lemma}
\label{lem2}
For $\beta$ large,
\begin{equation}
\lim_{n \rightarrow \infty} \mbox{e}^{-\beta n} \nu^+_\beta
 \left( \mbox{e}^{2\beta \sum_{i=-n}^n (1-\sigma(i))}
 \right)=0.
\end{equation} 
\end{lemma}   

{\bf Proof:}\\
Note that
\begin{equation}
\mbox{e}^{-\beta n} \nu^+_\beta \left(
\mbox{e}^{2\beta \sum_{i=-n}^n (1-\sigma(i))}
\right)
\leq \mbox{e}^{-\frac{\beta}{3} n} \nu^+_\beta \left(
 \mbox{e}^{\frac{11}{6} \beta \sum_{i=-n}^n(1-\sigma(i))} \right).
\end{equation}
Obviously, $\nu^+_\beta(f) = \mbox{lim}_\Lambda \ \mu^+_{\beta,\Lambda}(f)$
and we can represent $\mu^+_{\beta,\Lambda}(\mbox{e}^{c\beta \sum_{i=-n}^n
(1-\sigma(i))})$
for $\Lambda$ large, in terms of a contour representation. In particular,
using the cluster expansion we will prove that 
for $0<c<2$ and $\beta>\frac{\beta_0}{2-c}$
(for some large $\beta_0$),
$\limsup_n \frac{1}{n}
\log \nu^+_\beta \left( \mbox{e}^{c\beta \sum_{i=-n}^n (1-\sigma(i))}\right)$
is uniformly bounded. 
We adopt the notation of \cite{Sim} and refer to sections V.7 and V.8 therein for 
details. Let $\Gamma_0= \Gamma_0(\Lambda)$ denote the set of all Ising contours 
(i.e. contours on the dual lattice) corresponding to configurations in
volume $\Lambda$.
The dependence on $\Lambda$ will be understood and not explicitly kept
as all arguments will turn out to be uniform in the box $\Lambda$.
Using a contourrepresentation for the partition function, we can write
\begin{eqnarray}
\label{cluster}
 \nu^+_\beta \left(
\mbox{e}^{c\beta \sum(1-\sigma(i))}
\right) 
& \leq &
 \lim_\Lambda \ 
 \frac{\sum_\Gamma \prod_{\gamma \in \Gamma} \mbox{e}^{-2\beta |\gamma|}
 \mbox{e}^{2c\beta |int\gamma \cap [-n,n]|} }
 { \sum_\Gamma \prod_{\gamma \in \Gamma} \mbox{e}^{-2\beta |\gamma|} } 
   \\
& \leq &
\lim_\Lambda \
 \frac{\sum_\Gamma \prod_{\gamma \in \Gamma} \mbox{e}^{-(2-c)\beta |\gamma|}
    }
 { \sum_\Gamma \prod_{\gamma \in \Gamma} \mbox{e}^{-2\beta |\gamma|} }   
\end{eqnarray}
where the sums run over all families of mutually  disjoint 
contours ($\Gamma = \{ \gamma_1, \ldots, \gamma_k: \ \gamma_i \in \Gamma_0 
\ \mbox{and} \ \gamma_i \cap \gamma_j = \emptyset, \ \forall i,j=1,
\ldots,k\}$), 
$int\gamma$ denotes the set of sites in the interior of $\gamma$ 
and we have used that $|int \gamma \cap [-n,\ldots,n]|$ $\leq \frac{1}{2}|\gamma|$.

We will now take the logarithm of the sums appearing in the RHS of 
(\ref{cluster}).
The cluster expansion enables us to write this logarithm again as 
a sum but now over  connected families of contours 
where every
contour can appear more than once. We therefore change the notation (still
following \cite{Sim}) and go over from sets of contours ($\Gamma$) 
to multi-indices $A$.

We call a map
$ A: \Gamma_0 \rightarrow \N$
a multi-index. $A(\gamma)$ has to be interpreted as the number of times that the
contour $\gamma$ appears. Next we 
define (for $0 \leq r < 2 $)
\begin{equation}
z_r(\gamma):=\mbox{e}^{-(2-r)\beta |\gamma|},
\end{equation}

\begin{equation}
z_r^A \ :=\ \prod_{\gamma \in \Gamma_0} z_r(\gamma)^{A(\gamma)}, 
\end{equation} 
then
\begin{equation}
\label{clust}
\log \sum_\Gamma \prod_{\gamma \in \Gamma} z_r^{|\gamma|} =
\sum_A a^T(A) \ z_r^A
\end{equation}
for suitable coefficients $a^T(A)$ (see \cite{Sim} p466 for their exact 
expression).
Using  (\ref{clust}) and (\ref{cluster}) we obtain
\begin{equation}
\log \ \nu^+_\beta \left(
 \mbox{e}^{c\beta \sum(1-\sigma(i))} 
 \right)
\leq
\lim_\Lambda \ 
\sum_{A} a^T(A) \left[
  z_c^A - z_0^A \right]
\end{equation}
and the multi-indices in the sum must 
give non-zero weight to at least one contour that intersects $[-n,n]$ 
otherwise the expression in the square brackets is zero.
Using translation invariance, taking the limit $\Lambda \uparrow \Z^2$,
deviding by $2n+1$ and taking the limit 
$n \rightarrow \infty$ we get
\begin{equation}
\limsup_{n \rightarrow \infty} \frac{1}{2n+1} \log \nu^+_\beta \left(
\mbox{e}^{c\beta \sum_{i=-n}^n (1-\sigma(i))} \right) 
\leq 
2\ \sum_{A \ni 0} a^T (A) z_c^A.
\end{equation} 
This sum converges and is of order 
$\mbox{e}^{-(2-c)\beta}$ for $c<2$ and $\beta>\frac{\beta_0}{2-c}$.
In other words, for $\beta$ sufficiently large we have that
\begin{equation}
\nu^+_\beta \left( \mbox{e}^{ \frac{11}{6}\beta \sum_{i=-n}^n (1-\sigma(i))}
\right) = O\left(\mbox{e}^{n {\rm e}^{- \frac{1}{6}\beta}} \right)
\end{equation}
and hence
\begin{equation}
\lim_{n \rightarrow \infty} \mbox{e}^{-\frac{1}{3}\beta n} \nu^+_\beta 
\left( \mbox{e}^{\frac{11}{6}\beta \sum_{i=-n}^n (1-\sigma(i))} \right) =0.
\end{equation}
\QED

\subsection{Open Problem}
Schonmann proved in \cite{Sch} that 
\begin{equation}
s(\nu^-_\beta|\nu^+_\beta) > 0.
\end{equation}
This implies that the following two assertions cannot be true
together: 
\begin{itemize}
\item
both phases $\nu^+_\beta$ and $\nu^-_\beta$
are consistent with $\Gamma$ 
\item
\begin{equation}
0=s(\nu^+_\beta |U) = \lim_V \frac{1}{|V|} 
S_V (\nu^+_\beta|\nu^-_\beta\circ\gamma_V ).
\end{equation} 
If the first assertion is false, then $\nu^+_\beta$ and $\nu^-_\beta$
are not almost Gibbsian (i.e. the set of continuity points of $\Gamma$
has $\nu^+_\beta$ and $\nu^-_\beta$ measure zero, see \cite{FP}). It is
however believed (though not proved) 
that $\nu^+_\beta$ and $\nu^-_\beta$ are actually almost
Gibbsian. In that case the second assertion must be false and the limit
$\lim_V \frac{1}{|V|} S_V (\nu^+_\beta |\gamma_V (\cdot|\omega ) )$ is not
the same for all boundary conditions $\omega$.
\end{itemize}

\end{document}